\documentclass{amsart}

\newtheorem{theorem}{Theorem}[section]
\newtheorem{lemma}[theorem]{Lemma}
\newtheorem{proposition}[theorem]{Proposition}
\newtheorem{corollary}[theorem]{Corollary}

\theoremstyle{definition}

\newtheorem{remark}[theorem]{Remark}

\usepackage{amscd,amssymb}

\begin{document}

\title[Defining Relations of $3\times 3$ Trace Algebra]
{Defining Relations of Minimal Degree of the Trace Algebra of
$3 \times 3$ Matrices}
\author[Francesca Benanti and Vesselin Drensky]
{Francesca Benanti and Vesselin Drensky}
\address{Dipartimento di Matematica ed Applicazioni, Universit\`a di Palermo,
Via Archirafi 34, 90123 Palermo, Italy}
\email{fbenanti@math.unipa.it}
\address{Institute of Mathematics and Informatics,
Bulgarian Academy of Sciences,
          1113 Sofia, Bulgaria}
\email{drensky@math.bas.bg}
\thanks{The research of the first author was partially supported by MIUR, Italy.}
\thanks{The research of the second author was partially supported by Grant
MI-1503/2005 of the Bulgarian National Science Fund.}

\begin{abstract}
The trace algebra $C_{nd}$ over a field of characteristic 0
is generated by all traces of products of
$d$ generic $n\times n$ matrices, $n,d\geq 2$.
Minimal sets of generators of $C_{nd}$ are known for $n=2$ and
$n=3$ for any $d$ as well as for $n=4$ and $n=5$ and $d=2$. The defining relations between
the generators are found for $n=2$ and any $d$ and for $n=3$, $d=2$ only.
Starting with the generating set of $C_{3d}$
given by Abeasis and Pittaluga in 1989, we have shown that the
minimal degree of the set of defining relations of $C_{3d}$
is equal to 7 for any $d\geq 3$.
We have determined all relations of minimal degree.
For $d=3$ we have also found
the defining relations of degree 8. The proofs are
based on methods of representation theory of the general linear
group and easy computer calculations with standard functions of
Maple.
\end{abstract}

\maketitle

\section*{Introduction}

Let $K$ be any field of characteristic 0. All vector spaces, tensor products,
algebras considered in this paper are over $K$. Let
$X_i=\left(x_{pq}^{(i)}\right)$, $p,q=1,\ldots,n$, $i=1,\ldots,d$,
be $d$ generic $n\times n$ matrices. We consider the
pure (or commutative) trace algebra $C_{nd}$
generated by all traces of products
$\text{\rm tr}(X_{i_1}\cdots X_{i_k})$. The algebra
$C_{nd}$ coincides with the algebra of invariants of the general
linear group $GL_n=GL_n(K)$ acting by simultaneous conjugation on
$d$ matrices of size $n\times n$. General
results of invariant theory of classical groups imply that the
algebra $C_{nd}$ is finitely generated.
Theory of PI-algebras provides upper bounds for the
generating sets of the algebras $C_{nd}$.
The Nagata-Higman theorem states that the polynomial
identity $x^n=0$ implies the identity $x_1\cdots x_N=0$ for some
$N=N(n)$. If $N$ is minimal with this property, then $C_{nd}$ is
generated by traces of products
$\text{\rm tr}(X_{i_1}\cdots X_{i_k})$
of degree $k\leq N$.
This estimate is sharp if $d$ is sufficiently
large. A description of the defining relations of $C_{nd}$ is
given by the Razmyslov-Procesi theory \cite{R, P} in the language
of ideals of the group algebras of symmetric groups. For a
background on the algebras of matrix invariants
see, e.g. \cite{F, DF} and for computation aspects of the theory see
\cite{D2}.

Explicit minimal sets of generators of $C_{nd}$ and
the defining relations between them are found in few cases only.
It is well known that, in the Nagata-Higman theorem,
$N(2)=3$, $N(3)=6$, and $N(4)=10$,
which gives bounds for the degrees of the generators
of the algebras $C_{2d}$, $C_{3d}$, and $C_{4d}$, respectively.
Nevertheless, the defining relations of $C_{nd}$ are explicitly given
for $n=2$ and any $d$, see e.g. \cite{DF} for details, and for $n=3$, $d=2$,
see the comments below. For $n=3$, $d\geq 3$ and $n\geq 4$ and $d\geq 2$, nothing is known about the
concrete form of the defining relations with respect to fixed minimal systems of generators.

Teranishi \cite{T1} found a system of 11 generators of $C_{32}$.
It follows from his description that, with respect to these generators,
$C_{32}$ can be defined by a single relation of degree 12.
The explicit (but very complicated) form of
the relation was found by Nakamoto \cite{N}, over $\mathbb Z$,
with respect to a slightly different system of generators.
Abeasis and Pittaluga \cite{AP} found
a system of generators of $C_{3d}$, for any $d\geq 2$, in terms of
representation theory of the symmetric and general linear groups,
in the spirit of its usage in theory of PI-algebras.
Aslaksen, Drensky and Sadikova \cite{ADS}
gave the defining relation of $C_{32}$ with respect to the set
found in \cite{AP}. Their relation is much simpler
than that in \cite{N}. For $C_{42}$, a set of generators was found by Teranishi \cite{T1, T2}
and a minimal set by Drensky and Sadikova \cite{DS}, in terms of the approach in \cite{AP}.
Djokovi\'{c} \cite{D} gave another minimal set of 32 generators
of $C_{42}$ consisting of traces of products only. He found also a minimal set
of 173 generators of $C_{52}$.

As usually in invariant theory, the determination of
generators and defining relations is simpler, if one has some additional information about
the algebras of invariants. In particular, it is very useful to know the Hilbert (or Poincar\'e)
series of the algebra. Again, the picture is completely clear for $n=2$. The only other cases,
when the Hilbert series are explicitly given, are $n=3$, $d=2$ (Teranishi \cite{T1}) and $d=3$
(Berele and Stembridge \cite{BS}), $n=4$, $d=2$
(Teranishi \cite{T2} (with some typos) and corrected by
Berele and Stembridge \cite{BS}). Recently Djokovi\'{c} \cite{D}
has calculated also the Hilbert series of $C_{52}$ and $C_{62}$.

The minimal generating set of $C_{3d}$ given in \cite{AP} consists of
\[
g=g(d)=\frac{1}{240}d(5d^5+19d^4-5d^3+65d^2+636)
\]
homogeneous trace polynomials $u_1,\ldots,u_g$ of degree $\leq 6$. In more detail,
the number of polynomials $u_i$ of degree $k$ is $g_k$, and
\[
g_1=d,g_2=\frac{1}{2}(d+1)d,g_3=\frac{1}{3}d(d^2+2),
g_4=\frac{1}{24}(d+1)d(d-1)(5d-6),
\]
\[
g_5=\frac{1}{30}d(d-1)(d-2)(3d^2+4d+6),
g_6=\frac{1}{48}(d+2)(d+1)d(d-1)(d^2-3d+4).
\]
Hence $C_{3d}$ is isomorphic to the factor algebra $K[y_1,\ldots,y_g]/I$.
Defining $\text{\rm deg}(y_i)=\text{\rm deg}(u_i)$,
the ideal $I$ is homogeneous. For $d=3$, surprisingly,
the comparison of the Hilbert series of $C_{33}\cong K[y_1,\ldots,y_g]/I$ given in \cite{BS},
with the Hilbert series of $K[y_1,\ldots,y_g]$, gives that
any homogeneous minimal system of generators of the ideal $I$ contains
no elements of degree $\leq 6$, three elements of degree 7 and 30 elements of degree 8.
The purpose of the present paper is to find the defining relations
of minimal degree for $C_{3d}$ and any $d\geq 3$, with respect to
the generating set in \cite{AP}. It has turned out
that the minimal degree of the relations is equal to 7
for all $d\geq 3$, and there are a lot of relations of degree 7.
(Compare with the single relation of degree 12
in the case $d=2$.) The dimension of the vector space
of relations of degree 7 is equal to
\[
r_7=r_7(d)=\frac{2}{7!}(d+1)d(d-1)(d-2)(41d^3-86d^2+114d-360).
\]
For $d=3$ we have computed also the homogeneous relations of degree 8.
The defining relations are given in the language of representation theory of $GL_d$.
There is a simple algorithm which gives the explicit form of all relations of degree 7,
and of degree 8 for $d=3$.
The proofs involve basic representation theory of $GL_d$ and develop further ideas
of \cite{ADS, DS} and our recent paper \cite{BD}
combined with computer calculations with Maple. In the case
$d=3$ we have used essentially the Hilbert
series of $C_{33}$ from \cite{BS}
which has allowed to reduce the number of computations.
Our methods are quite general and we believe that they can be successfully
used for further investigation of generic trace algebras
and other algebras close to them.

\section{Preliminaries}

In what follows, we fix $n=3$ and $d\geq 3$ and denote by $X_1,\ldots,X_d$
the $d$ generic $3\times 3$ matrices. Often, when the value of $d$ is clear
from the context, we shall denote $C_{3d}$ by $C$.
It is a standard trick to
replace the generic matrices with generic traceless
matrices. We express $X_i$ in the form
\[
X_i=\frac{1}{3}\text{\rm tr}(X_i)e+x_i,\quad i=1,\ldots,d,
\]
where $e$ is the identity $3\times 3$ matrix and $x_i$ is a
generic traceless matrix.
Then
\begin{equation}\label{replacing with traceless matrices}
C_{3d}\cong K[\text{\rm tr}(X_1),\ldots,\text{\rm tr}(X_d)]\otimes C_0,
\end{equation}
where the algebra $C_0$ is generated by the traces of products
$\text{\rm tr}(x_{i_1}\cdots x_{i_k})$, $k\leq 6$. Hence the problem for
the defining relations of $C$ can be replaced by a similar problem for $C_0$.

As in the case of ``ordinary'' generic matrices, without loss of generality we may
replace $x_1$ by a generic traceless diagonal matrix. Changing the variables
$x_{pp}^{(i)}$, we may assume that
\begin{equation}\label{first matrix}
x_1=\left( \begin{array}{ccc}
                x_{11}^{(1)} & 0 & 0 \\
                0 & x_{22}^{(1)} & 0 \\
                0 & 0 & -(x_{11}^{(1)}+x_{22}^{(1)})\\
             \end{array} \right),
\end{equation}
\begin{equation}\label{other matrices}
x_i=\left( \begin{array}{ccc}
                x_{11}^{(i)} & x_{12}^{(i)} & x_{13}^{(i)} \\
                x_{21}^{(i)} & x_{22}^{(i)} & x_{23}^{(i)} \\
                x_{31}^{(i)} & x_{32}^{(i)} & -(x_{11}^{(i)}+x_{22}^{(i)})\\
             \end{array} \right), \quad i=2,\ldots,d.
\end{equation}
Till the end of the paper we fix these $d$
generic traceless matrices. Let $C_0^+=\omega(C_0)$ be the augmentation ideal
of $C_0$. It consists of all trace polynomials $f(x_1,\ldots,x_d)\in C_0$ without constant terms,
i.e., satisfying the condition
$f(0,\ldots,0)=0$. Any minimal system of generators of $C_0$ lying in
$C_0^+$ forms a basis of the vector space $C_0^+$ modulo $(C_0^+)^2$.
Conversely, if a system of polynomials $f_1,\ldots,f_g$
forms a basis of $C_0^+$ modulo $(C_0^+)^2$, and each $f_i$ is a linear combination of
traces of products $\text{\rm tr}(x_{i_1}\cdots x_{i_k})$
then it is a minimal generating system of $C_0$. The algebra $C=C_{3d}$
is $\mathbb Z$-graded assuming that the trace $\text{\rm tr}(X_{i_1}\cdots X_{i_k})$ is of degree $k$,
and this grading is inherited by $C_0$. Similarly, $C$ (and also $C_0$)
has a more precise ${\mathbb Z}^d$-multigrading induced by the condition
that $X_1,\ldots,X_d$ are, respectively, of multidegree
$(1,0,\ldots,0,0),\ldots,(0,0,\ldots,0,1)$. The considerations below, stated for
the $\mathbb Z$-grading hold also for the ${\mathbb Z}^d$-multigrading.
The numbers $g_1,g_2,\ldots,g_6$ of elements of degree $1,2,\ldots,6$, respectively,
in any homogeneous minimal system of generators is an invariant of $C$.
Any homogeneous minimal system $\{f_1,\ldots,f_h\}$ of generators of $C_0$
consists of $g_2,\ldots,g_6$ elements of degree $2,\ldots,6$, and $h=g_2+\cdots+g_6$. Hence
\[
C_0\cong K[z_1,\ldots,z_h]/J,
\]
with isomorphism defined by $z_j+J\to f_j$, $j=1,\ldots,h$.
If $u_j(z_1,\ldots,z_h)$, $j=1,\ldots,r$, is a system of generators of the ideal $J$,
then $u_j(f_1,\ldots,f_h)=0$, $j=1,\ldots,r$, is a system of defining relations of $C_0$
with respect to the system of generators
$\{f_1,\ldots,f_h\}$. Any homogeneous system of polynomials
in $J$ (where, by definition $\text{\rm deg}(z_j)=\text{\rm deg}(f_j)$)
which forms a basis of the vector space $J$ modulo the subspace $JK[z_1,\ldots,z_h]^+$,
is a minimal system of generators of the ideal $J$. We denote by $r_k$
the number of elements of degree $k$ in such a system. Clearly,
$r_k$ is the dimension of the homogeneous component of degree $k$ of
the vector space $J/JK[z_1,\ldots,z_h]^+$.

Now we summarize the necessary background on representation theory
of $GL_d$. We refer e.g. to \cite{M} for general facts
and to \cite{D1}  for applications in the spirit of the problems
considered here.
All $GL_d$-modules which appear in this paper are
completely reducible and are direct sums of irreducible polynomial
modules.
The irreducible polynomial representations of $GL_d$ are
indexed by partitions $\lambda=(\lambda_1,\ldots,\lambda_d)$,
$\lambda_1\geq \cdots\geq \lambda_d\geq 0$. We denote by
$W(\lambda)=W_d(\lambda)$ the
corresponding irreducible $GL_d$-module, assuming that $W_d(\lambda)=0$
if $\lambda_{d+1}\not=0$. The group $GL_d$ acts
in the natural way on the $d$-dimensional vector space
$K\cdot x_1+\cdots+K\cdot x_d$ and this action is extended diagonally on
the free associative algebra $K\langle x_1,\ldots,x_d\rangle$.

The module $W(\lambda)\subset K\langle x_1,\ldots,x_d\rangle$ is generated by a
unique, up to a multiplicative constant,
homogeneous element $w_{\lambda}$ of degree $\lambda_j$
with respect to $x_j$,
called the highest weight vector of $W(\lambda)$.
It is characterized by the following property.

\begin{lemma}\label{criterion for hwv}
Let $1\leq i<j\leq d$ and let
$\Delta_{ij}$ be the derivation of $K\langle x_1,\ldots,x_d\rangle$
defined by $\Delta_{ij}(x_j)=x_i$, $\Delta_{ij}(x_k)=0$, $k\not=j$.
If $w(x_1,\ldots,x_d) \in K\langle x_1,\ldots,x_d\rangle$
is multihomogeneous of degree $\lambda=(\lambda_1,\ldots,\lambda_d)$,
then $w(x_1,\ldots,x_d)$ is a highest weight vector for some
$W(\lambda)$ if and only if $\Delta_{ij}(w(x_1,\ldots,x_d))=0$
for all $i<j$. Equivalently, $w(x_1,\ldots,x_d)$ is a highest weight vector
for $W(\lambda)$ if and only if
\[
g_{ij}(w(x_1,\ldots,x_d))=w(x_1,\ldots,x_d),\quad 1\leq i<j\leq d,
\]
where $g_{ij}$ is the linear operator of the $d$-dimensional vector space which
sends $x_j$ to $x_i+x_j$ and fixes the other $x_k$.
\end{lemma}

\begin{proof} The lemma is a partial case of a result by
De Concini, Eisenbud, and Procesi \cite{DEP}, see also
Almkvist, Dicks, and Formanek \cite{ADF}.
In the version which we need, the first part of the lemma was established
by Koshlukov \cite{K}. The equivalence follows from the fact that
the kernel of any locally nilpotent derivation $\Delta$ coincides with the fixed points
of the related exponential automorphism
$\exp(\Delta)=1+\Delta/1!+\Delta^2/2!+\cdots$,
and $g_{ij}=\exp(\Delta_{ij})$.
\end{proof}

If $W_i$, $i=1,\ldots,m$,
are $m$ isomorphic copies of the $GL_d$-module $W(\lambda)$
and $w_i\in W_i$ are highest weight
vectors, then the highest weight vector of any submodule $W(\lambda)$
of the direct sum $W_1\oplus\cdots\oplus W_m$ has the form
$\xi_1w_1+\cdots+\xi_mw_m$ for some $\xi_i\in K$.
Any $m$ linearly independent highest weight vectors can serve
as a set of generators of the $GL_d$-module $W_1\oplus\cdots\oplus W_m$.

It is convenient to work with an explicit copy of $W(\lambda)$ in $K\langle x_1,\ldots,x_d\rangle$
obtained in the following way. Let
\[
s_k(x_1,\ldots,x_k)=\sum_{\sigma\in S_k}\text{\rm sign}(\sigma)x_{\sigma(1)}\cdots x_{\sigma(k)}
\]
be the standard polynomial of degree $k$.
(Clearly,
\[
s_2(x_1,x_2)=x_1x_2-x_2x_1=[x_1,x_2]
\]
is the commutator of $x_1$ and $x_2$.)
If the lengths of the columns of the diagram of $\lambda$
are, respectively, $k_1,\ldots,k_p$, $p=\lambda_1$, then
\begin{equation}\label{canonical hwv}
w_{\lambda}=w_{\lambda}(x_1,\ldots,x_{k_1})=s_{k_1}(x_1,\ldots,x_{k_1})\cdots s_{k_p}(x_1,\ldots,x_{k_p})
\end{equation}
is the highest weight vector of a submodule $W(\lambda)\subset K\langle x_1,\ldots,x_d\rangle$.
Sometimes we shall write $w_{\lambda}=w_{\lambda}(x_1,\ldots,x_d)$, even when $k_1<d$.

Recall that the $\lambda$-tableau
\[
T=(a_{ij}),\quad a_{ij}\in \{1,\ldots,d\},\quad i=1,\ldots,d,\quad j=1,\ldots,\lambda_i,
\]
is semistandard if its entries do not decrease from left to right in rows
and increase from top to bottom in columns. The following lemma gives a basis of the vector subspace
$W(\lambda)\subset K\langle x_1,\ldots,x_d\rangle$. It also provides
an algorithm to construct this basis.

\begin{lemma}\label{basis of module}
Let $T=(a_{ij})$ be a semistandard $\lambda$-tableau such that its $i$-th row
contains $b_{i,i}$ times $i$, $b_{i,i+1}$ times $i+1$, $\ldots$, $b_{i,d}$ times $d$.
Let $w(x_1,\ldots,x_d)$ be the highest weight vector of
$W(\lambda)\subset K\langle x_1,\ldots,x_d\rangle$ and let
\[
u_T(x_{11},x_{12},\ldots,x_{1d},x_{22},\ldots,x_{2d},\ldots,x_{dd})
\]
be the multihomogeneous component of degree $b_{iq}$ in $x_{iq}$, $q=i,i+1,\ldots,d$, of the polynomial
\[
w(x_{11}+x_{12}+\cdots+x_{1d},x_{22}+\cdots+x_{2d},\ldots,x_{dd}).
\]
When $T$ runs on the set of semistandard $\lambda$-tableaux, the polynomials
\[
v_T=v_T(x_1,\ldots,x_d)=u_T(x_1,x_2,\ldots,x_d,x_2,\ldots,x_d,\ldots,x_d)
\]
form a basis of the vector space $W(\lambda)$.
\end{lemma}

\begin{proof}
By standard Vandermonde arguments, the polynomial
$u_T(x_{11},x_{12},\ldots,x_{dd})$ is a linear combination of some
$w(\sum\alpha_{1j}x_{1j},\ldots,\sum\alpha_{dj}x_{dj})$, $\alpha_{ij}\in K$.
Hence $v_T(x_1,\ldots,x_d)$ is a linear combination of
$w(\sum\alpha_{1j}x_j,\ldots,\sum\alpha_{dj}x_j)$ and belongs to the $GL_d$-module $W(\lambda)$
generated by $w(x_1,\ldots,x_d)$. Without loss of generality,
it is sufficient to consider the case when $W(\lambda)$ is generated by the element
(\ref{canonical hwv}). The polynomial
$w(x_{11}+x_{12}+\cdots+x_{1d},x_{22}+\cdots+x_{2d},\ldots,x_{dd})$ is a product of evaluations
\[
s_{k_j}(x_{11}+x_{12}+\cdots+x_{1d},x_{22}+\cdots+x_{2d},\ldots,x_{k_j,k_j}+\cdots+x_{k_j,d})
\]
of standard polynomials. Hence $u_T(x_{11},x_{12},\ldots,x_{dd})$
is a linear combination of monomials starting with
$x_{\sigma(1),q_1}\cdots x_{\sigma(k_1),q_{k_1}}$.
We order the variables by $x_1>\cdots>x_d$ and consider the lexicographic
order of $K\langle x_1,\ldots,x_d\rangle$. The first column of the semistandard tableau $T$
contains $a_{11}<\cdots<a_{k1}$ and in each row $a_{i1}\leq a_{ij}$ for all $j=2,\ldots,\lambda_i$.
This easily implies that the leading monomial of $v_T$ starts with $x_{a_{11}}x_{a_{21}}\cdots x_{a_{k_11}}$.
Hence the first $k_1$ variables in the leading monomial are indexed by the entries of the first column
of the tableau. Continuing in the same way, we obtain that the next $k_2$ variables in the leading monomial
are indexed with the entries of the second row, etc. Hence, for fixed $\lambda$, the leading monomial of
$v_T$ determines completely the tableau $T$, and the polynomials $v_T$ are linearly independent.
Since $W(\lambda)$ has a basis which is in 1-1 correspondence with the semistandard $\lambda$-tableaux,
and the number of $v_T$ is equal to the number of semistandard tableaux, we obtain that
the polynomials $v_T$ form a basis
of $W(\lambda)$.
\end{proof}

If $W$ is a $GL_d$-submodule or a factor module of
$K\langle x_1,\ldots,x_d\rangle$, then $W$ inherits the ${\mathbb Z}^d$-grading of
$K\langle x_1,\ldots,x_d\rangle$.
Recall that the Hilbert series of $W$ with respect to
its ${\mathbb Z}^d$-multigrading
is defined as the formal power series
\[
H(W,t_1,\ldots,t_d)=\sum_{k_i\geq 0}
\dim(W^{(k_1,\ldots,k_d)})t_1^{k_1}\cdots t_d^{k_d},
\]
with coefficients equal to the dimensions of the homogeneous
components $W^{(k_1,\ldots,k_d)}$ of degree $(k_1,\ldots,k_d)$.
It plays the role of the $GL_d$-character of $W$: If
\[
W\cong \sum_{\lambda}m(\lambda)W(\lambda),
\]
then
\[
H(W,t_1,\ldots,t_d)=\sum_{\lambda}m(\lambda)S_{\lambda}(t_1,\ldots,t_d),
\]
where $S_{\lambda}=S_{\lambda}(t_1,\ldots,t_d)$
is the Schur function associated with
$\lambda$, and the multiplicities $m(\lambda)$ are determined
by $H(W,t_1,\ldots,t_d)$.
One of the possible ways to introduce Schur functions is via
Vandermonde-like determinants. For a partition
$\mu=(\mu_1,\ldots,\mu_d)$, define the determinant
\[
V(\mu_1,\ldots,\mu_d)=\left\vert\begin{matrix}
t_1^{\mu_1}&t_2^{\mu_1}&\cdots&t_d^{\mu_1}\\
t_1^{\mu_2}&t_2^{\mu_2}&\cdots&t_d^{\mu_2}\\
\vdots&\vdots&\ddots&\vdots\\
t_1^{\mu_d}&t_2^{\mu_d}&\cdots&t_d^{\mu_d}\\
\end{matrix}\right\vert.
\]
Then the Schur function is
\[
S_{\lambda}(t_1,\ldots,t_d)=
\frac{V(\lambda_1+d-1,\lambda_2+d-2,\ldots,\lambda_{d-1}+1,\lambda_d)}
{V(d-1,d-2,\ldots,1,0)}.
\]
The dimension of $W(\lambda)$ is given by the formula
\begin{equation}\label{dim of W}
\text{\rm dim}(W_d(\lambda))=\prod_{1\leq i<j\leq d}\frac{\lambda_i-\lambda_j+j-i}{j-i}.
\end{equation}
The decomposition of the tensor product $W_d(\lambda)\otimes W_d(\mu)$
of two irreducible $GL_d$-modules $W_d(\lambda)$ and $W_d(\mu)$ is given by
the Littlewood-Richardson rule. We shall need it in one case only:
\begin{equation}\label{usage of LRrule}
W_3(2^2)\otimes W_3(2^2)\cong W_3(4^2)\oplus W_3(4,3,1)\oplus W_3(4,2^2).
\end{equation}
But even in this case we can check the equality (\ref{usage of LRrule})
verifying directly the equality of symmetric functions
\[
H(W_3(2^2)\otimes W_3(2^2),t_1,t_2,t_3)=S_{(2^2)}^2(t_1,t_2,t_3)
\]
\[
=S_{(4^2)}(t_1,t_2,t_3)
+S_{(4,3,1)}(t_1,t_2,t_3)
+S_{(4,2^2)}(t_1,t_2,t_3).
\]
In all other cases it will be sufficient to use the Young rule which
is a partial case of the Littlewood-Richardson one:
\begin{equation}\label{Young rule A}
W_d(\lambda_1,\ldots,\lambda_d)\otimes W_d(p)\cong
\sum W_d(\lambda_1+p_1,\ldots,\lambda_d+p_d),
\end{equation}
where the sum runs on all nonnegative integers $p_1,\ldots,p_d$
such that $p_1+\cdots+p_d=p$ and
$\lambda_i\geq \lambda_{i+1}+p_{i+1}$, $i=1,\ldots,d-1$,
and its dual version
\begin{equation}\label{Young rule B}
W_d(\lambda_1,\ldots,\lambda_d)\otimes W_d(1^p)\cong
\sum W_d(\lambda_1+\varepsilon_1,\ldots,\lambda_d+\varepsilon_d),
\end{equation}
where the sum is on all partitions
$(\lambda_1+\varepsilon_1,\ldots,\lambda_d+\varepsilon_d)$
such that $\varepsilon_i=0,1$, $\varepsilon_1+\cdots+\varepsilon_d=p$.

In the $q$-th symmetric tensor power
\[
W^{\otimes_sq}=\underbrace{W\otimes_s\cdots\otimes_sW}_{q\text{ \rm times}}
\]
of the $GL_d$-module $W$,
we identify the tensors $w_{\sigma(1)}\otimes\cdots\otimes w_{\sigma(q)}$
and $w_1\otimes\cdots\otimes w_q$, $\sigma\in S_q$.
If $W=W_1\oplus\cdots \oplus W_k$, then
\begin{equation}\label{symmetric powers}
W^{\otimes_sq}
=\bigoplus W_1^{\otimes_sq_1}\otimes\cdots\otimes W_k^{\otimes_sq_k},\quad q_1+\cdots+q_k=q.
\end{equation}
There is no general
combinatorial rule for the decomposition of the symmetric tensor powers
of $W_d(\lambda)$. We shall need the following partial results
due to Thrall \cite{Th}, see also \cite{M}:
\begin{equation}\label{first Thrall}
K[W_d(2)]=\sum_{q\geq 0}W_d(2)^{\otimes_sq}
=\sum W_d(2\lambda_1,2\lambda_2,\ldots,2\lambda_d),
\end{equation}
where the sum is over all partitions $\lambda$,
\begin{equation}\label{second Thrall}
W_d(p)\otimes_sW_d(p)=\sum_{0\leq k\leq p/2}W_d(p+2k,p-2k),
\end{equation}
\begin{equation}\label{third Thrall}
W_d(1^p)\otimes_sW_d(1^p)=\sum_{0\leq k\leq p/2}W_d(2^{p-2k},1^{4k}).
\end{equation}
Besides, we shall need the decomposition
\begin{equation}\label{symmetric square of W22}
W_3(2^2)\otimes_sW(2^2)=W_3(4^2)\oplus W_3(4,2^2).
\end{equation}
The easiest way to check (\ref{symmetric square of W22}) is
to use (\ref{usage of LRrule}), hence
\[
W_3(2^2)\otimes_sW(2^2)\subset W_3(4^2)\oplus W_3(4,3,1)\oplus W_3(4,2^2).
\]
Therefore
\[
W_3(2^2)\otimes_sW(2^2)= \varepsilon_1W_3(4^2)\oplus
\varepsilon_2W_3(4,3,1)\oplus \varepsilon_3W_3(4,2^2),
\]
$\varepsilon_i=0,1$. The Hilbert series of $W_3(2^2)\otimes_sW_3(2^2)$
contains the summand $t_1^4t_2^4$ and $S_{(4^2)}$
is the only Schur function among
$S_{(4^2)},S_{(4,3,1)},S_{(4,2^2)}$ which contains $t_1^4t_2^4$.
This implies that $W_3(4^2)$ participates in the decomposition of
$W_3(2^2)\otimes_sW_3(2^2)$ and $\varepsilon_1=1$.
Finally, we apply dimension arguments:
\[
\text{\rm dim}(W_3(2^2)\otimes_sW(2^2))=\text{\rm dim}(W_3(4^2))+
\varepsilon_2\text{\rm dim}(W_3(4,3,1))+\varepsilon_3\text{\rm dim}(W_3(4,2^2)).
\]
Since
\[
\text{\rm dim}(W_3(2^2))=6,\quad
\text{\rm dim}(W_3(2^2)\otimes_sW(2^2))=\binom{6+1}{2}=21,
\]
\[
\text{\rm dim}(W_3(4^2))=\text{\rm dim}(W_3(4,3,1))=15,\quad
\text{\rm dim}(W_3(4,2^2))=6,
\]
we obtain the only possibility $\varepsilon_2=0$, $\varepsilon_3=1$.

The action of $GL_d$ on $K\langle x_1,\ldots,x_d\rangle$
is inherited by the algebras $C_{3d}$ and $C_0$.
Now we discuss the approach of Abeasis and Pittaluga \cite{AP} for the special case $n=3$.
(Pay attention that the partitions in \cite{AP} are given in ``Francophone'' way, i.e., transposed to ours.)
The algebra $C_{3d}$ has a system of generators of degree $\leq 6$. Without
loss of generality we may assume that this system consists of traces of products
$\text{\rm tr}(X_{i_1}\cdots X_{i_k})$.
Let $U_k$ be the subalgebra of $C_{3d}$ generated by all traces
$\text{\rm tr}(X_{i_1}\cdots X_{i_l})$ of degree $l\leq k$.
Clearly, $U_k$ is also a $GL_d$-submodule of $C_{3d}$.
Let $C_{3d}^{(k+1)}$ be the homogeneous component of degree $k+1$
of $C_{3d}$. Then the intersection $U_k\cap C_{3d}^{(k+1)}$
is a $GL_d$-module and has a complement
$G_{k+1}$ in $C_{3d}^{(k+1)}$,
which is the $GL_d$-module of the ``new'' generators of degree $k+1$.
We may assume that $G_{k+1}$ is a submodule of the $GL_d$-module
spanned by traces of products $\text{\rm tr}(X_{i_1}\cdots X_{i_{k+1}})$ of degree $k+1$.
The $GL_d$-module of the generators of $C_{3d}$ is
\[
G=G_1\oplus G_2\oplus \cdots\oplus G_6.
\]

\begin{proposition}\label{generating module of C}
{\rm (Abeasis and Pittaluga \cite{AP})}
The $GL_d$-module $G$ of the generators of $C_{3d}$ decomposes as
\[
G=W(1)\oplus W(2)\oplus W(3)\oplus W(1^3)\oplus W(2^2)\oplus
W(2,1^2)
\]
\[
\oplus W(3,1^2)\oplus W(2^2,1)\oplus W(1^5)\oplus
W(3^2)\oplus W(3,1^3).
\]
Each module $W(\lambda)\subset G$ is generated by the ``canonical'' highest weight
vector $\text{\rm tr}(w_{\lambda}(X_1,\ldots,X_d))$, where $w_{\lambda}$
is given in {\rm (\ref{canonical hwv})}.
\end{proposition}

\begin{corollary}\label{dimensions of G}
The numbers $g_k$ of generators of degree $k\leq 6$
in any homogeneous minimal system of generators of $C_{3d}$ are
\[
g_1=d,g_2=\frac{1}{2}(d+1)d,g_3=\frac{1}{3}d(d^2+2),
g_4=\frac{1}{24}(d+1)d(d-1)(5d-6),
\]
\[
g_5=\frac{1}{30}d(d-1)(d-2)(3d^2+4d+6),
g_6=\frac{1}{48}(d+2)(d+1)d(d-1)(d^2-3d+4).
\]
The total number of generators is
\[
g=g(d)=\frac{1}{240}d(5d^5+19d^4-5d^3+65d^2+636)
\]
\end{corollary}

\begin{proof}
The number $g_k$ is equal to the dimension of the $GL_d$-submodule $G_k$ of
the $GL_d$-module $G$ of generators of $C_{3d}$. Applying the formula
(\ref{dim of W}) we obtain that the dimensions of the $GL_d$-modules
\[
W(1),W(2),W(3),W(1^3),W(2^2),W(2,1^2),
\]
\[
W(3,1^2),W(2^2,1),W(1^5),W(3^2),W(3,1^3).
\]
are, respectively,
\[
d,\binom{d+1}{2},\binom{d+2}{3},\binom{d}{3},
\frac{d}{2}\binom{d+1}{3},3\binom{d+1}{4},
\]
\[
6\binom{d+2}{5},d\binom{d+1}{4},
\binom{d}{5},
3\binom{d+2}{4}\binom{d+1}{2},10\binom{d+2}{6}.
\]
This easily implies the results, because
\[
\text{\rm dim}(G_3)=\text{\rm dim}(W(3))+\text{\rm dim}(W(1^3)),
\text{\rm dim}(G_4)=\text{\rm dim}(W(2^2))+\text{\rm dim}(W(2,1^2)),
\]
\[
\text{\rm dim}(G_5)=\text{\rm dim}(W(3,1^2))+\text{\rm dim}(W(2^2,1))
+\text{\rm dim}(W(1^5)),
\]
\[
\text{\rm dim}(G_6)=\text{\rm dim}(W(3^2))+\text{\rm dim}(W(3,1^3)),
\]
and $g=g_1+g_2+\cdots+g_6$.
\end{proof}

In the sequel we shall need the Hilbert series of
$C_{33}$ calculated by Berele and Stembridge \cite{BS}:
\[
H(C_{33},t_1,t_2,t_3)= \frac{p(t_1,t_2,t_3)}{q(t_1,t_2,t_3)},
\]
where
\[
p=1-e_2+e_3+e_1e_3+e_2^2+e_1^2e_3-e_2e_3-2e_1e_2e_3+e_3^2+e_2^2e_3
\]
\[
-e_1^2e_2e_3+2e_1^2e_3^2+e_1^3e_3^2+e_2^2e_3^2-e_1^2e_2e_3^2-e_1e_3^3-2e_1e_2^2e_3^2
\]
\[
+2e_2e_2e_3^3-e_2^3e_3^2+e_1^3e_3^3+2e_1^2e_2e_3^3-2e_1e_3^4-e_1^2e_3^4+e_1e_2^2e_3^3
\]
\[
+e_2e_3^4-e_2^3e_3^3-2e_2^2e_3^5+e_1e_2^2e_3^4+2e_1e_2e_3^5-e_3^6-e_2^2e_3^5
\]
\[
+e_1e_3^6-e_2e_3^6-e_1^2e_3^6-e_3^7+e_1e_3^7-e_3^8,
\]
\[
q=\!\!\left(\prod_{i=1}^3(1-t_i)(1-t_i^2)(1-t_i^3)\right)\!\!\!
\left(\prod_{1\leq i<j\leq
3}(1-t_it_j)^2(1-t_i^2t_j)(1-t_it_j^2)\right)\!\! (1-t_1t_2t_3),
\]
and
\[
e_1=t_1+t_2+t_3,\quad e_2=t_1t_2+t_1t_3+t_2t_3,\quad e_3=t_1t_2t_3
\]
are the elementary symmetric polynomials in three variables.
Since the Hilbert series of the tensor product is equal to the product
of the Hilbert series of the factors, and
\[
H(K[\text{\rm \rm tr}(X_1),\ldots,\text{\rm \rm tr}(X_d)],t_1,\ldots,t_d)=
\frac{1}{(1-t_1)\cdots(1-t_d)},
\]
(\ref{replacing with traceless matrices}) implies that
\[
H(C_{33},t_1,t_2)=\frac{H(C_0,t_1,t_2,t_3)}{(1-t_1)(1-t_2)(1-t_3)}.
\]
In this way, for $d=3$,
\begin{equation}\label{Hilbert series of C0}
H(C_0,t_1,t_2,t_3)=
(1-t_1)(1-t_2)(1-t_3)H(C_{33},t_1,t_2).
\end{equation}

\section{The symmetric algebra of the generators}

We consider the symmetric algebra
\[
S=K[G_2\oplus \cdots\oplus G_6]
\]
of the $GL_d$-module of the generators of the algebra $C_0$. Clearly,
the grading and the $GL_d$-module structure of
$G_2\oplus \cdots\oplus G_6$ induce a grading and the structure of a
$GL_d$-module also on $S$. The defining relations of the algebra $C_0$
are in the square of the augmentation ideal $\omega(S)$ of $S$.
Since we are interested in the defining relations of degree 7 for $C_0$
for any $d\geq 3$ and of degree 8 for $d=3$,
we shall decompose the homogeneous components of degree 7, respectively, 8 of
the ideal $\omega^2(S)$ into a sum of irreducible $GL_d$-, respectively, $GL_3$-modules.
Then we shall find explicit generators of those irreducible components
which may give rise to relations.

\begin{lemma}\label{multiplicities of degree 7}
The following $GL_d$-module isomorphisms hold:
\begin{equation}\label{W2xW2}
W(2)\otimes_sW(2)\cong W(4)\oplus W(2^2),
\end{equation}
\begin{equation}\label{W31xW2}
W(3)\otimes W(2)\cong W(5)\oplus W(4,1)\oplus W(3,2),
\end{equation}
\begin{equation}\label{W32xW2}
W(1^3)\otimes W(2)\cong W(3,1^2)\oplus W(2,1^3),
\end{equation}
\begin{equation}\label{W41xW2}
W(2^2)\otimes W(2)\cong W(4,2)\oplus W(3,2,1)\oplus W(2^3),
\end{equation}
\begin{equation}\label{W42xW2}
W(2,1^2)\otimes W(2)\cong W(4,1^2)\oplus W(3,2,1)\oplus W(3,1^3)
\oplus W(2^2,1^2),
\end{equation}
\begin{equation}\label{W31xW31}
W(3)\otimes_sW(3)\cong W(6)\oplus W(4,2),
\end{equation}
\begin{equation}\label{W31xW32}
W(3)\otimes W(1^3)\cong W(4,1^2)\oplus W(3,1^3),
\end{equation}
\begin{equation}\label{W32xW32}
W(1^3)\otimes_sW(1^3)\cong W(2^3)\oplus W(2,1^4),
\end{equation}
\begin{equation}\label{W2xW2xW2}
W(2)\otimes_sW(2)\otimes_sW(2)\cong W(6)\oplus W(4,2)\oplus W(2^3),
\end{equation}
\begin{equation}\label{W51xW2}
W(3,1^2)\otimes W(2)\cong W(5,1^2)\oplus W(4,2,1)\oplus W(4,1^3)
\oplus W(3^2,1)\oplus W(3,2,1^2),
\end{equation}
\begin{equation}\label{W52xW2}
W(2^2,1)\otimes W(2)\cong W(4,2,1)\oplus W(3,2^2)\oplus W(3,2,1^2)
\oplus W(2^3,1),
\end{equation}
\begin{equation}\label{W53xW2}
W(1^5)\otimes W(2)\cong W(3,1^4)\oplus W(2,1^5),
\end{equation}
\begin{equation}\label{W41xW31}
W(2^2)\otimes W(3)\cong W(5,2)\oplus W(4,2,1)\oplus W(3,2^2),
\end{equation}
\begin{equation}\label{W41xW32}
W(2^2)\otimes W(1^3)\cong W(3^2,1)\oplus W(3,2,1^2)\oplus W(2^2,1^3),
\end{equation}
\begin{equation}\label{W42xW31}
W(2,1^2)\otimes W(3)\cong W(5,1^2)\oplus W(4,2,1)\oplus W(4,1^3)
\oplus W(3,2,1^2),
\end{equation}
\begin{equation}\label{W42xW32}
W(2,1^2)\otimes W(1^3)\cong W(3,2^2)\oplus W(3,2,1^2)
\end{equation}
\begin{equation}\nonumber
\oplus W(3,1^4) \oplus W(2^3,1)\oplus W(2^2,1^3)\oplus W(2,1^5),
\end{equation}
\begin{equation}\label{W31xW2xW2}
W(3)\otimes(W(2)\otimes_sW(2))\cong W(7)\oplus W(6,1)
\end{equation}
\begin{equation}\nonumber
\oplus 2W(5,2) \oplus W(4,3)\oplus W(4,2,1)\oplus W(3,2^2),
\end{equation}
\begin{equation}\label{W32xW2xW2}
W(1^3)\otimes(W(2)\otimes_sW(2))\cong W(5,1^2)\oplus W(4,1^3)\oplus W(3^2,1)
\oplus W(2^2,1^3).
\end{equation}
\end{lemma}

\begin{proof}
Equations (\ref{W31xW31}), (\ref{W32xW32}) are partial cases of (\ref{first Thrall}), (\ref{second Thrall}),
respectively;
(\ref{W2xW2}) and (\ref{W2xW2xW2}) follow from (\ref{third Thrall}).
Equations (\ref{W31xW2}), (\ref{W32xW2}), (\ref{W41xW2}),
(\ref{W42xW2}), (\ref{W31xW32}), (\ref{W51xW2}), (\ref{W52xW2}),
(\ref{W53xW2}), (\ref{W41xW31}), and (\ref{W42xW31}) are obtained from
(\ref{Young rule A}); (\ref{W41xW32}) and (\ref{W42xW32}) follow from
(\ref{Young rule B}). Finally,
(\ref{W31xW2xW2}) and (\ref{W32xW2xW2}) are calculated applying,
respectively, (\ref{Young rule A}) and (\ref{Young rule B})
to (\ref{W2xW2}).
\end{proof}

\begin{proposition}\label{decomposition of low degree}
The homogeneous components $(\omega^2(S))^{(k)}$ of degree $k\leq 7$ of the square $\omega^2(S)$
of the augmentation ideal of the symmetric algebra of
$G_2\oplus\cdots\oplus G_6$ decomposes as
\[
(\omega^2(S))^{(4)}=W(4)\oplus W(2^2),
\]
\[
(\omega^2(S))^{(5)}=W(5)\oplus W(4,1)\oplus W(3,2)\oplus W(3,1^2)\oplus W(2,1^3),
\]
\[
(\omega^2(S))^{(6)}=2W(6)\oplus 3W(4,2)\oplus 2W(4,1^2)
\]
\[
\oplus 2W(3,2,1)
\oplus 2W(3,1^3)\oplus 3W(2^3)\oplus W(2^2,1^2)\oplus W(2,1^4),
\]
\[
(\omega^2(S))^{(7)}=W(7)\oplus W(6,1)\oplus 3W(5,2)\oplus 3W(5,1^2)
\]
\[
\oplus W(4,3)\oplus 5W(4,2,1)
\oplus 3W(4,1^3)\oplus 3W(3^2,1)\oplus 4W(3,2^2)
\]
\[
\oplus 6W(3,2,1^2)\oplus 2W(3,1^4)\oplus 2W(2^3,1)\oplus 3W(2^2,1^3)\oplus 2W(2,1^5).
\]
\end{proposition}

\begin{proof}
By (\ref{symmetric powers}),
the homogeneous component of degree $k$ of
$\omega^2(S)$ has the form
\[
(\omega^2(S))^{(k)}
=\bigoplus G_2^{\otimes_sq_2}\otimes\cdots\otimes G_6^{\otimes_sq_6},\quad 2q_2+\cdots+6q_6=k, q_2+\cdots+q_6\geq 2.
\]
The decomposition of $G_2,\ldots,G_6$ is given in
Proposition \ref{generating module of C}. Again, the equality
(\ref{symmetric powers}) implies that $(\omega^2(S))^{(k)}=0$ for $k<4$.
For $k=4$ we use (\ref{W2xW2}) and obtain
\[
(\omega^2(S))^{(4)}= G_2\otimes_sG_2=W(2)\otimes_sW(2)=W(4)\oplus W(2^2).
\]
For $k=5$, we derive
\[
(\omega^2(S))^{(5)}= G_3\otimes G_2=(W(3)\oplus W(1^3))\otimes W(2)
\]
and the decomposition follows from (\ref{W31xW2}) and (\ref{W32xW2}).
For $k=6$ we have
\[
(\omega^2(S))^{(6)}=G_4\otimes G_2\oplus G_3\otimes_sG_3
\oplus G_2\otimes_s G_2\otimes_sG_2
\]
and we use the decompositions given in
(\ref{W41xW2}) -- (\ref{W2xW2xW2}). The decomposition of
$(\omega^2(S))^{(7)}$ is obtained in a similar way and makes use of
(\ref{W51xW2}) -- (\ref{W32xW2xW2}).
\end{proof}

\begin{proposition}\label{hwv of degree up to 6}
The following elements of $S=K[G_2\oplus\cdots\oplus G_6]$
are highest weight vectors:

\noindent For $\lambda=(2,1^3)$:
\[
w=\sum_{\sigma\in S_4}\text{\rm sign}(\sigma)
\text{\rm tr}(s_3(x_{\sigma(1)},x_{\sigma(2)},x_{\sigma(3)}))
\text{\rm \rm tr}(x_{\sigma(4)}x_1);
\]
For $\lambda=(3,1^3)$:
\[
w_1=\sum_{\sigma\in S_4}\text{\rm sign}(\sigma)
\text{\rm tr}(s_3(x_{\sigma(1)},x_{\sigma(2)},x_{\sigma(3)})x_1
+s_3(x_1,x_{\sigma(1)},x_{\sigma(2)})x_{\sigma(3)})
\text{\rm tr}(x_{\sigma(4)}x_1),
\]
\[
w_2=\sum_{\sigma\in S_4}\text{\rm sign}(\sigma)
\text{\rm tr}(x_1^2x_{\sigma(1)})
\text{\rm tr}(s_3(x_{\sigma(2)},x_{\sigma(3)},x_{\sigma(4)})),
\]
For $\lambda=(2^2,1^2)$:
\[
w=\sum_{\sigma\in S_4,\tau\in S_2}\text{\rm sign}(\sigma\tau)
\text{\rm tr}(s_3(x_{\sigma(1)},x_{\sigma(2)},x_{\sigma(3)})x_{\tau(1)}
\]
\[
+s_3(x_{\tau(1)},x_{\sigma(1)},x_{\sigma(2)})x_{\sigma(3)})
\text{\rm tr}(x_{\sigma(4)}x_{\tau(2)}),
\]
For $\lambda=(2,1^4)$:
\[
w=\sum_{\sigma\in S_5}\text{\rm sign}(\sigma)
\text{\rm tr}(s_3(x_{\sigma(1)},x_{\sigma(2)},x_{\sigma(3)}))
\text{\rm tr}(s_3(x_{\sigma(4)},x_{\sigma(5)},x_1)).
\]
For $\lambda=(2,1^3), (2^2,1^2), (2,1^4)$, every highest weight vector
$w\in W(\lambda)\subset \omega^2(S)$ is equal, up to a multiplicative constant,
to the corresponding $w$. For $\lambda=(3,1^3)$
the highest weight vectors are linear combinations of $w_1$ and $w_2$.
\end{proposition}

\begin{proof}
By Proposition \ref{generating module of C}, the submodules
$W(2), W(3),W(1^3),W(2,1^2)$ of $G_2\oplus\cdots\oplus G_6$ are generated,
respectively, by
\[
v_1=\text{\rm tr}(x_1^2),v_2=\text{\rm tr}(x_1^3),v_3=\text{\rm tr}(s_3(x_1,x_2,x_3)),
v_4=\text{\rm tr}(s_3(x_1,x_2,x_3)x_1).
\]
The trace polynomials $\text{\rm tr}(x_1x_2)$ and
$\text{\rm tr}(s_3(x_1,x_2,x_3)x_4+s_3(x_2,x_3,x_4)x_1)$
are linearizations of $v_1$ and $v_4$, respectively.
Hence the trace polynomials $w,w_1,w_2$ defined in the proposition
belong to the ideal $\omega^2(S)$.
They have the necessary skew-symmetries, as the ``canonical'' highest weight
vectors $w_{\lambda}$ from (\ref{canonical hwv}). Hence all
$w$ and $w_1,w_2$ are highest weight vectors. Clearly, they are nonzero
in $S$ and their number coincides with the multiplicities of $W(\lambda)$
in the decomposition of $\omega^2(S)$. Hence, every highest weight vector
$w\in W(\lambda)\subset \omega^2(S)$ can be expressed as their linear combination.
\end{proof}

\begin{proposition}\label{hwv of degree 7}
The following elements of $S=K[G_2\oplus\cdots\oplus G_6]$
are highest weight vectors:

\noindent For $\lambda=(4,1^3)$:
\[
w_1=(\text{\rm tr}(s_3(x_1,x_2,x_3)(x_1x_4+x_4x_1))-\text{\rm tr}(s_3(x_1,x_2,x_4)(x_1x_3+x_3x_1))
\]
\[
+\text{\rm tr}(s_3(x_1,x_3,x_4)(x_1x_2+x_2x_1))
+3\text{\rm tr}(s_3(x_2,x_3,x_4)x_1^2))\text{\rm tr}(x_1^2)
\]
\[
+5(-\text{\rm tr}(s_3(x_1,x_2,x_3)x_1^2)\text{\rm tr}(x_1x_4)
\]
\[
+\text{\rm tr}(s_3(x_1,x_2,x_4)x_1^2)\text{\rm tr}(x_1x_3)
-\text{\rm tr}(s_3(x_1,x_3,x_4)x_1^2)\text{\rm tr}(x_1x_2)),
\]
\[
w_2=(\text{\rm tr}(s_3(x_1,x_2,x_3)x_4)-\text{\rm tr}(s_3(x_1,x_2,x_4)x_3)+\text{\rm tr}(s_3(x_1,x_3,x_4)x_2)
\]
\[
+3\text{\rm tr}(s_3(x_2,x_3,x_4)x_1))\text{\rm tr}(x_1^3)
+4(-\text{\rm tr}(s_3(x_1,x_2,x_3)x_1)\text{\rm tr}(x_1^2x_4)
\]
\[
+\text{\rm tr}(s_3(x_1,x_2,x_4)x_1)\text{\rm tr}(x_1^2x_3)
-\text{\rm tr}(s_3(x_1,x_3,x_4)x_1)\text{\rm tr}(x_1^2x_2)),
\]
\[
w_3=(\text{\rm tr}(s_3(x_2,x_3,x_4)))\text{\rm tr}(x_1^2)-\text{\rm tr}(s_3(x_1,x_3,x_4)))\text{\rm tr}(x_1x_2)
\]
\[
 +\text{\rm tr}(s_3(x_1,x_2,x_4)))\text{\rm tr}(x_1x_3)
 -\text{\rm tr}(s_3(x_1,x_2,x_3)))\text{\rm tr}(x_1x_4))\text{\rm tr}(x_1^2).
\]

\noindent For $\lambda=(3,2^2)$:
\[
w_1=\sum_{\sigma\in S_3}\text{\rm sign}(\sigma)
\text{\rm tr}(s_3(x_1,x_2,x_3)x_{\sigma(1)}x_{\sigma(2)})\text{\rm tr}(x_1x_{\sigma(3)}),
\]
\[
w_2=\text{\rm tr}(s_3(x_1,x_2,x_3)x_1)\text{\rm tr}(s_3(x_1,x_2,x_3)),
\]
\[
w_3=\text{\rm tr}([x_1,x_2]^2)\text{\rm tr}(x_1x_3^2)+\text{\rm tr}([x_1,x_3]^2)\text{\rm tr}(x_1x_2^2)
+\text{\rm tr}([x_2,x_3]^2)\text{\rm tr}(x_1^3)
\]
\[
-\text{\rm tr}([x_1,x_2][x_1,x_3])\text{\rm tr}(x_1(x_2x_3+x_3x_2))
\]
\[
+2\text{\rm tr}([x_1,x_2][x_2,x_3])\text{\rm tr}(x_1^2x_3)
-2\text{\rm tr}([x_1,x_3][x_2,x_3])\text{\rm tr}(x_1^2x_2),
\]
\[
w_4=\text{\rm tr}(x_1^3)(\text{\rm tr}(x_2^2)\text{\rm tr}(x_3^2)-\text{\rm tr}^2(x_2x_3))
+\text{\rm tr}(x_1x_2^2)(\text{\rm tr}(x_1^2)\text{\rm tr}(x_3^2)-\text{\rm tr}^2(x_1x_3))
\]
\[
+\text{\rm tr}(x_1x_3^2)(\text{\rm tr}(x_1^2)\text{\rm tr}(x_2^2)-\text{\rm tr}^2(x_1x_2))
\]
\[
+2\text{\rm tr}(x_1^2x_2)(-\text{\rm tr}(x_1x_2)\text{\rm tr}(x_3^2)+\text{\rm tr}(x_1x_3)\text{\rm tr}(x_2x_3))
\]
\[
+2\text{\rm tr}(x_1^2x_3)(-\text{\rm tr}(x_1x_3)\text{\rm tr}(x_2^2)+\text{\rm tr}(x_1x_2)\text{\rm tr}(x_2x_3))
\]
\[
+\text{\rm tr}(x_1(x_2x_3+x_3x_2))(-\text{\rm tr}(x_1^2)\text{\rm tr}(x_2x_3)
+\text{\rm tr}(x_1x_2)\text{\rm tr}(x_1x_3)).
\]

\noindent For $\lambda=(3,2,1^2)$:
\[
w_1=(\text{\rm tr}(s_3(x_1,x_2,x_3)(x_2x_4+x_4x_2)) -\text{\rm tr}((s_3(x_1,x_2,x_4)(x_2x_3+x_3x_2))
\]
\[
+4\text{\rm tr}(s_3(x_2,x_3,x_4)(x_1x_2+x_2x_1))+ 2\text{\rm tr}(s_3(x_1,x_3,x_4)x_2^2))\text{\rm tr}(x_1^2)
\]
\[
+(-\text{\rm tr}(s_3(x_1,x_2,x_3)(x_1x_4+x_4x_1)) +\text{\rm tr}(s_3(x_1,x_2,x_4)(x_1x_3+x_3x_1))
\]
\[
-6\text{\rm tr}(s_3(x_1,x_3,x_4)(x_1x_2+x_2x_1)) -8\text{\rm tr}(s_3(x_2,x_3,x_4)x_1^2))\text{\rm tr}(x_1x_2)
\]
\[
+5(-\text{\rm tr}(s_3(x_1,x_2,x_3)(x_1x_2+x_2x_1))\text{\rm tr}(x_1x_4)
\]
\[
+\text{\rm tr}(s_3(x_1,x_2,x_4)(x_1x_2+x_2x_1))\text{\rm tr}(x_1x_3))
\]
\[
+10(\text{\rm tr}(s_3(x_1,x_2,x_3)x_1^2)\text{\rm tr}(x_2x_4)
-\text{\rm tr}(s_3(x_1,x_2,x_4)x_1^2)\text{\rm tr}(x_2x_3)
\]
\[
+\text{\rm tr}(s_3(x_1,x_3,x_4)x_1^2)\text{\rm tr}(x_2^2)),
\]
\[
w_2=\text{\rm tr}(x_1^2)(-\text{\rm tr}(s_3(x_1,x_2,x_3)[x_2,x_4])
-\text{\rm tr}(s_3(x_2,x_3,x_4)[x_1,x_2])
\]
\[
+\text{\rm tr}(s_3(x_1,x_2,x_4)[x_2,x_3])
+ 3\text{\rm tr}(s_3(x_2,x_3,x_4)[x_1,x_2]))
\]
\[
+\text{\rm tr}(x_1x_2)(\text{\rm tr}(s_3(x_1,x_2,x_3)[x_1,x_4])
- \text{\rm tr}(s_3(x_1,x_3,x_4)[x_1,x_2])
\]
\[
-\text{\rm tr}(s_3(x_1,x_2,x_4)[x_1,x_3])-\text{\rm tr}(s_3(x_1,x_3,x_4)[x_1,x_2]))
\]
\[
+ 3\text{\rm tr}(x_1x_3)\text{\rm tr}(s_3(x_1,x_2,x_4)[x_1,x_2])
- 3\text{\rm tr}(x_1x_4)\text{\rm tr}(s_3(x_1,x_2,x_3)[x_1,x_2]),
\]
\[
w_3=\text{\rm tr}([x_1,x_2]^2)\text{\rm tr}(s_3(x_1,x_3,x_4))
- \text{\rm tr}([x_1,x_2][x_1,x_3])\text{\rm tr}(s_3(x_1,x_2,x_4))
\]
\[
+ \text{\rm tr}([x_1,x_2][x_1,x_4])\text{\rm tr}(s_3(x_1,x_2,x_3)),
\]
\[
w_4=-2\text{\rm tr}(s_3(x_2,x_3,x_4)x_2)\text{\rm tr}(x_1^3)
+ 2(\text{\rm tr}(s_3(x_1,x_3,x_4)x_2)
\]
\[
+ \text{\rm tr}(s_3(x_2,x_3,x_4)x_1))\text{\rm tr}(x_1^2x_2)
-  2\text{\rm tr}(s_3(x_1,x_2,x_4)x_2)\text{\rm tr}(x_1^2x_3)
\]
\[
+2\text{\rm tr}(s_3(x_1,x_2,x_3)x_2)\text{\rm tr}(x_1^2x_4)
- 2\text{\rm tr}(s_3(x_1,x_3,x_4)x_1)\text{\rm tr}(x_1x_2^2)
\]
\[
+ \text{\rm tr}(s_3(x_1,x_2,x_4)x_1)\text{\rm tr}(x_1(x_2x_3+x_3x_2))
- \text{\rm tr}(s_3(x_1,x_2,x_3)x_1)\text{\rm tr}(x_1(x_2x_4+x_4x_2)),
\]
\[
w_5=\text{\rm tr}(s_3(x_1,x_2,x_3)x_1)\text{\rm tr}(s_3(x_1,x_2,x_4))
- \text{\rm tr}(s_3(x_1,x_2,x_4)x_1)\text{\rm tr}(s_3(x_1,x_2,x_3)),
\]
\[
w_6=(\text{\rm tr}(x_1^2)\text{\rm tr}(x_2^2)-\text{\rm tr}(x_1x_2)^2)\text{\rm tr}(s_3(x_1,x_3,x_4))
\]
\[
+ (-\text{\rm tr}(x_1^2)\text{\rm tr}(x_2x_3)
+\text{\rm tr}(x_1x_2)\text{\rm tr}(x_1x_3))\text{\rm tr}(s_3(x_1,x_2,x_4))
\]
\[
+ (\text{\rm tr}(x_1^2)\text{\rm tr}(x_2x_4)
-\text{\rm tr}(x_1x_2)\text{\rm tr}(x_1x_4))\text{\rm tr}(s_3(x_1,x_2,x_3)).
\]

\noindent For $\lambda=(3,1^4)$:
\[
w_1=\text{\rm tr}(s_5(x_1,x_2,x_3,x_4,x_5))\text{\rm tr}(x_1^2),
\]
\[
w_2=\sum_{\sigma\in S_5}\text{\rm sign}(\sigma)
\text{\rm tr}(s_3(x_1,x_{\sigma(2)},x_{\sigma(3)})x_1)\text{\rm tr}(s_3(x_1,x_{\sigma(4)},x_{\sigma(5)})),
\quad \sigma(1)=1.
\]

\noindent For $\lambda=(2^3,1)$:
\[
w_1=\text{\rm tr}(s_3(x_2,x_3,x_4)[x_2,x_3])\text{\rm tr}(x_1^2)
\]
\[
-(\text{\rm tr}(s_3(x_1,x_3,x_4)[x_2,x_3]) +\text{\rm tr}(s_3(x_2,x_3,x_4)[x_1,x_3]))\text{\rm tr}(x_1x_2)
\]
\[+(\text{\rm tr}(s_3(x_1,x_2,x_4)[x_2,x_3]) +\text{\rm tr}(s_3(x_2,x_3,x_4)[x_1,x_2]))\text{\rm tr}(x_1x_3)
\]
\[
-\text{\rm tr}(s_3(x_1,x_2,x_3)[x_2,x_3])\text{\rm tr}(x_1x_4)
+\text{\rm tr}(s_3(x_1,x_3,x_4)[x_1,x_3])\text{\rm tr}(x_2^2)
\]
\[ -(\text{\rm tr}(s_3(x_1,x_2,x_4)[x_1,x_3]) + \text{\rm tr}(s_3(x_1,x_3,x_4)[x_1,x_2]))\text{\rm tr}(x_2x_3)
\]
\[+\text{\rm tr}(s_3(x_1,x_2,x_3)[x_1,x_3])\text{\rm tr}(x_2x_4)
+\text{\rm tr}(s_3(x_1,x_2,x_4)[x_1,x_2])\text{\rm tr}(x_3^2)
\]
\[
-\text{\rm tr}(s_3(x_1,x_2,x_3)[x_1,x_2])\text{\rm tr}(x_3x_4),
\]
\[
w_2=(-3(\text{\rm tr}(s_3(x_1,x_2,x_3)x_4) + \text{\rm tr}(s_3(x_2,x_3,x_4)x_1))+\text{\rm tr}(s_3(x_1,x_3,x_4)x_2)
\]
\[+ \text{\rm tr}(s_3(x_2,x_3,x_4)x_1))-\text{\rm tr}(s_3(x_1,x_2,x_4)x_3)
\]
\[
+ \text{\rm tr}(s_3(x_2,x_3,x_4)x_1)))\text{\rm tr}(s_3(x_1,x_2,x_3))
\]
\[
+4\text{\rm tr}(s_3(x_1,x_2,x_3)x_3)\text{\rm tr}(s_3(x_1,x_2,x_4))
\]
\[
-4\text{\rm tr}(s_3(x_1,x_2,x_3)x_2)\text{\rm tr}(s_3(x_1,x_3,x_4))
\]
\[
+4\text{\rm tr}(s_3(x_1,x_2,x_3)x_1)\text{\rm tr}(s_3(x_2,x_3,x_4)).
\]

\noindent For $\lambda=(2^2,1^3)$:
\[
w_1=(\text{\rm tr}([x_1,x_4][x_2,x_5]) + \text{\rm tr}([x_1,x_5][x_2,x_4])
\]
\[
+2\text{\rm tr}([x_1,x_2][x_4,x_5]) - 2\text{\rm tr}([x_1,x_5][x_2,x_4]))\text{\rm tr}(s_3(x_1,x_2,x_3))
\]
\[+ (-\text{\rm tr}([x_1,x_3][x_2,x_5]) - \text{\rm tr}([x_1,x_5][x_2,x_3])-2\text{\rm tr}([x_1,x_2][x_3,x_5])
\]
\[+2\text{\rm tr}([x_1,x_5][x_2,x_3]))\text{\rm tr}(s_3(x_1,x_2,x_4))+ (\text{\rm tr}([x_1,x_3][x_2,x_4]) +
\]
\[\text{\rm tr}([x_1,x_4][x_2,x_3])+2\text{\rm tr}([x_1,x_2][x_3,x_4]) -
\]
\[2\text{\rm tr}([x_1,x_4][x_2,x_3]))\text{\rm tr}(s_3(x_1,x_2,x_5))
\]
\[
+ 3\text{\rm tr}([x_1,x_2][x_2,x_5])\text{\rm tr}(s_3(x_1,x_3,x_4))
- 3\text{\rm tr}([x_1,x_2][x_2,x_4])\text{\rm tr}(s_3(x_1,x_3,x_5))
\]
\[
+ 3\text{\rm tr}([x_1,x_2][x_2,x_3])\text{\rm tr}(s_3(x_1,x_4,x_5))
- 3\text{\rm tr}([x_1,x_2][x_1,x_5])\text{\rm tr}(s_3(x_2,x_3,x_4))
\]
\[
+ 3\text{\rm tr}([x_1,x_2][x_1,x_4])\text{\rm tr}(s_3(x_2,x_3,x_5))
- 3\text{\rm tr}([x_1,x_2][x_1,x_3])\text{\rm tr}(s_3(x_2,x_4,x_5))
\]
\[
+ 3\text{\rm tr}([x_1,x_2]^2)\text{\rm tr}(s_3(x_3,x_4,x_5)),
\]
\[
w_2=\sum_{\sigma\in S_5}\sum_{\tau\in S_2}\text{\rm sign}(\sigma\tau)
\text{\rm tr}(s_3(x_{\tau(1)},x_{\sigma(1)},x_{\sigma(2)}))
(\text{\rm tr}(s_3(x_{\sigma(3)},x_{\sigma(4)},x_{\sigma(5)})x_{\tau(2)})\]
\[
+\sum_{\sigma\in S_5}\sum_{\tau\in S_2}\text{\rm sign}(\sigma\tau)
\text{\rm tr}(s_3(x_{\tau(1)},x_{\sigma(1)},x_{\sigma(2)}))
(\text{\rm tr}(s_3(x_{\tau(2)},x_{\sigma(3)},x_{\sigma(4)})x_{\sigma(5)}),
\]
\[
w_3=\text{\rm tr}(s_3(x_1,x_2,x_3))(\text{\rm tr}(x_1x_4)\text{\rm tr}(x_2x_5)
-\text{\rm tr}(x_1x_5)\text{\rm tr}(x_2x_4))
\]
\[
+ \text{\rm tr}(s_3(x_1,x_2,x_4))(-\text{\rm tr}(x_1x_3)\text{\rm tr}(x_2x_5)
+\text{\rm tr}(x_1x_5)\text{\rm tr}(x_2x_3))
\]
\[
+ \text{\rm tr}(s_3(x_1,x_2,x_5))(\text{\rm tr}(x_1x_3)\text{\rm tr}(x_2x_4)
-\text{\rm tr}(x_1x_4)\text{\rm tr}(x_2x_3))
\]
\[
+ \text{\rm tr}(s_3(x_1,x_3,x_4))(\text{\rm tr}(x_1x_2)\text{\rm tr}(x_2x_5)-\text{\rm tr}(x_1x_5)\text{\rm tr}(x_2^2))
\]
\[
+ \text{\rm tr}(s_3(x_1,x_3,x_5))(-\text{\rm tr}(x_1x_2)\text{\rm tr}(x_2x_4)
+ \text{\rm tr}(x_1x_4)\text{\rm tr}(x_2^2))
\]
\[
+ \text{\rm tr}(s_3(x_1,x_4,x_5))(\text{\rm tr}(x_1x_2)\text{\rm tr}(x_2x_3)
-\text{\rm tr}(x_1x_3)\text{\rm tr}(x_2^2))
\]
\[
+ \text{\rm tr}(s_3(x_2,x_3,x_4))(\text{\rm tr}(x_1x_2)\text{\rm tr}(x_1x_5)
-\text{\rm tr}(x_2x_5)\text{\rm tr}(x_1^2))
\]
\[
+ \text{\rm tr}(s_3(x_2,x_3,x_5))(-\text{\rm tr}(x_1x_2)\text{\rm tr}(x_1x_4)
+ \text{\rm tr}(x_2x_4)\text{\rm tr}(x_1^2))
\]
\[
+ \text{\rm tr}(s_3(x_2,x_4,x_5))(\text{\rm tr}(x_1x_2)\text{\rm tr}(x_1x_3)
-\text{\rm tr}(x_2x_3)\text{\rm tr}(x_1^2))
\]
\[
+ \text{\rm tr}(s_3(x_3,x_4,x_5))(\text{\rm tr}(x_1^2)\text{\rm tr}(x_2^2)-\text{\rm tr}(x_1x_2)^2).
\]

\noindent For $\lambda=(2,1^5)$:
\[
w_1=\sum_{\sigma\in S_6}\text{\rm sign}(\sigma)
\text{\rm tr}(s_5(x_{\sigma(1)},x_{\sigma(2)},x_{\sigma(3)},x_{\sigma(4)},x_{\sigma(5)}))
\text{\rm tr}(x_1x_{\sigma(6)}),
\]
\[
w_2=\sum_{\sigma\in S_6}\text{\rm sign}(\sigma)
\text{\rm tr}(s_3(x_{\sigma(1)},x_{\sigma(2)},x_{\sigma(3)})x_1))
\text{\rm tr}(s_3(x_{\sigma(4)},x_{\sigma(5)},x_{\sigma(6)}))
\]
\[
+\sum_{\sigma\in S_6}\text{\rm sign}(\sigma)
\text{\rm tr}(s_3(x_1,x_{\sigma(1)},x_{\sigma(2)})x_{\sigma(3)}))
\text{\rm tr}(s_3(x_{\sigma(4)},x_{\sigma(5)},x_{\sigma(6)})).
\]
For $\lambda=(4,1^3), (3,2^2), (3,2,1^2), (3,1^4), (2^3,1), (2^2,1^3),(2,1^5)$,
every highest weight vector $w\in W(\lambda)\subset \omega^2(S)$ is equal
to a linear combination of $w_i$.
\end{proposition}

\begin{proof} The computations are similar to those in the proof of
Proposition \ref{hwv of degree up to 6}. For a fixed $\lambda\vdash 7$,
the number of $w_i$ can be calculated from the equations
(\ref{W51xW2}) -- (\ref{W32xW2xW2}) and
Proposition \ref{decomposition of low degree}. The concrete form of $w_i$
is found using Lemmas \ref{criterion for hwv} and \ref{basis of module}.
We shall demonstrate the process in one case only.

Let $\lambda=(3,2^2)$. Consider the tensor product
$W(3)\otimes(W(2)\otimes_sW(2))\subset\omega^2(S)$. The equation
(\ref{W31xW2xW2}) gives that the module $W(3,2^2)$ participates with
multiplicity 1. More precisely, applying (\ref{Young rule A}) to
(\ref{W2xW2}), we see that $W(3,2^2)$ appears as a submodule of the
component $W(3)\otimes W(2^2)$ of $W(3)\otimes(W(2)\otimes_sW(2))$.
Since $W(3)\subset G_3$ and $W(2)=G_2$ are generated by $\text{\rm
tr}(x_1^3)$ and $\text{\rm tr}(x_1^2)$, respectively, Lemma
\ref{basis of module} gives that they have bases
\[
\{\text{\rm tr}(x_i^3),\text{\rm tr}(x_i^2x_j), i\not=j,\quad
\text{\rm tr}(x_i(x_jx_k+x_kx_j)), i<j<k\},
\]
\[
\{\text{\rm tr}(x_i^2), \text{\rm tr}(x_ix_j),\quad i\not=j\}.
\]
The submodule $W(2^2)$ of $W(2)\otimes_sW(2)$ is generated by
\[
u(x_1,x_2)=\text{\rm tr}(x_1^2)\text{\rm tr}(x_2^2)-\text{\rm tr}^2(x_1x_2).
\]
(Direct verification shows that $u(x_1,x_2+x_1)=u(x_1,x_2)$
and we apply Lemma \ref{criterion for hwv}.)
Its partial linearization in $x_2$ is
\[
v(x_1,x_2,x_3)=\text{\rm tr}(x_1^2)\text{\rm tr}(x_2x_3)
-\text{\rm tr}(x_1x_2)\text{\rm tr}(x_1x_3).
\]
We are looking for an element $w=w(x_1,x_2,x_3)$ in
$W(3)\otimes W(2^2)\subset W(3)\otimes(W(2)\otimes_sW(2))$ which is homogeneous
of multidegree $(3,2^2)$ and satisfies the conditions
$\Delta_{21}(w)=\Delta_{31}(w)=\Delta_{32}(w)=0$, where $\Delta_{ij}$
is the derivation from Lemma \ref{criterion for hwv}.
All such elements are of the form
\[
w=\zeta_1\text{\rm tr}(x_1^3)u(x_2,x_3)
+ \zeta_2\text{\rm tr}(x_1^2x_2)v(x_3,x_1,x_2)
+ \zeta_3\text{\rm tr}(x_1^2x_3)v(x_2,x_1,x_3)
\]
\[
+ \zeta_4\text{\rm tr}(x_1x_2^2)u(x_1,x_3)
+\zeta_5\text{\rm tr}(x_1(x_2x_3+x_3x_2))v(x_1,x_2,x_3)
+\zeta_6\text{\rm tr}(x_1x_3^2)u(x_1,x_2).
\]
Direct verifications show that
\[
\Delta_{21}(w)=(2\zeta_1+\zeta_2)\text{\rm tr}(x_1^3)v(x_3,x_1,x_2)
\]
\[
+(\zeta_2+2\zeta_4)\text{\rm tr}(x_1^2x_2)u(x_1,x_3)
+(-\zeta_3+2\zeta_5)\text{\rm tr}(x_1^2x_3)v(x_1,x_2,x_3),
\]
\[
\Delta_{31}(w)=(2\zeta_1+\zeta_3)\text{\rm tr}(x_1^3)v(x_2,x_1,x_3)
\]
\[
+(-\zeta_2+2\zeta_5)\text{\rm tr}(x_1^2x_2)v(x_1,x_2,x_3)
+(\zeta_3+2\zeta_6)\text{\rm tr}(x_1^2x_3)u(x_1,x_2),
\]
\[
\Delta_{32}(w)=(-\zeta_2+\zeta_3)\text{\rm tr}(x_1^2x_2)v(x_2,x_1,x_3)
\]
\[
+2(\zeta_4+\zeta_5)\text{\rm tr}(x_1x_2^2)v(x_1,x_2,x_3)
+(\zeta_5+\zeta_6)\text{\rm tr}(x_1(x_2x_3+x_3x_2))u(x_1,x_2).
\]
Hence we obtain the homogeneous linear system
\[
2\zeta_1+\zeta_2=\zeta_2+2\zeta_4=-\zeta_3+2\zeta_5=0,
\]
\[
2\zeta_1+\zeta_3=-\zeta_2+2\zeta_5=\zeta_3+2\zeta_6=0,
\]
\[
-\zeta_2+\zeta_3=+2(\zeta_4+\zeta_5)=\zeta_5+\zeta_6=0.
\]
Up to a multiplicative constant, the only solution of the system is
\[
\zeta_1=\zeta_4=\zeta_6=1,\quad \zeta_2=\zeta_3=-2,\quad\zeta_5=-1,
\]
which is equal to $w_4$.
In practice, in most of the cases we have used a slightly different algorithm to
determine $w_i$. We have considered $w_i$ with unknown coefficients.
Then we have evaluated it in the trace algebra $C_{3d}$ instead of in the symmetric algebra $S$,
in order to use the programs which we already had. Requiring that
\[
g_{kl}w(x_1,\dots,x_d)=w(x_1,\ldots,x_d), \quad 1\leq k<l\leq d,
\]
we have obtained the possible candidates for $w_i$. Since the number
of candidates has coincided with the number predicted in Proposition \ref{decomposition of low degree},
we have concluded that the $w_i$'s really are the needed highest weight vectors.
\end{proof}

\begin{proposition}\label{decomposition of degree 8}
For $d=3$, the homogeneous component $(\omega^2(S))^{(8)}$ of degree $8$ of the square $\omega^2(S)$
of the augmentation ideal of the symmetric algebra of
$G_2\oplus\cdots\oplus G_6$ decomposes as
\[
(\omega^2(S))^{(8)}=2W_3(8)\oplus W_3(7,1)\oplus 4W_3(6,2)\oplus 3W_3(6,1^2)
\]
\[
\oplus 2W_3(5,3)\oplus 6W_3(5,2,1)\oplus 4W_3(4^2)\oplus 7W_3(4,3,1)\oplus 9W_3(4,2^2)\oplus 4W_3(3^2,2).
\]
\end{proposition}

\begin{proof}
The considerations are similar to those in the proof of
Proposition \ref{decomposition of low degree} and involves the equtions
\[
W_3(3^2)\otimes W_3(2)\cong W_3(5,3)\oplus W_3(4,3,1)\oplus W_3(3^2,2),
\]
\[
W_3(3,1^2)\otimes W_3(3)\cong W_3(6,1^2)\oplus W_3(5,2,1)\oplus W_3(4,3,1),
\]
\[
W_3(3,1^2)\otimes W_3(1^3)\cong W_3(4,2^2),
\]
\[
W_3(2^2,1)\otimes W_3(3)\cong W_3(5,2,1)\oplus W_3(4,2^2),
\]
\[
W_3(2^2,1)\otimes W_3(1^3)\cong W_3(3^2,2),
\]
\[
W_3(2^2)\otimes_sW_3(2^2)\cong W_3(4^2)\oplus W_3(4,2^2),
\]
\[
W_3(2^2)\otimes W_3(2,1^2)\cong W_3(4,3,1)\oplus W_3(3^2,2),
\]
\[
W_3(2,1^2)\otimes_sW_3(2,1^2)\cong W_3(4,2^2),
\]
\[
W_3(2^2)\otimes(W_3(2)\otimes_sW_3(2))\cong W_3(6,2)
\]
\[
\oplus W_3(5,2,1)\oplus W_3(4^2)\oplus W_3(4,3,1)\oplus 2W_3(4,2^2),
\]
\[
W_3(2,1^2)\otimes(W_3(2)\otimes_sW_3(2))\cong W_3(6,1^2)\oplus W_3(5,2,1)\oplus W_3(4,3,1)\oplus W_3(3^2,2),
\]
\[
(W_3(3)\otimes_sW_3(3))\otimes W_3(2)\cong W_3(8)
\]
\[
\oplus W_3(7,1)\oplus 2W_3(6,2)\oplus W_3(5,3)\oplus W_3(5,2,1)
\]
\[
\oplus W_3(4^2)\oplus W_3(4,3,1)\oplus W_3(4,2^2),
\]
\[
(W_3(3)\otimes W_3(1^3))\otimes W_3(2)
\cong W_3(6,1^2)\oplus W_3(5,2,1) \oplus W_3(4,3,1),
\]
\[
(W_3(1^3)\otimes_sW_3(1^3))\otimes W_3(2)\cong W_3(4,2^2),
\]
\[
W_3(2)\otimes_sW_3(2)\otimes_sW_3(2)\otimes_sW_3(2)\cong
W_3(8)\oplus W_3(6,2)\oplus W_3(4^2)\oplus W_3(4,2^2).
\]
The proof of these equations uses the Young rule (\ref{Young rule A}) and (\ref{Young rule B}),
the formulas (\ref{first Thrall}) -- (\ref{third Thrall}), (\ref{usage of LRrule}),
and (\ref{symmetric square of W22}). For example,
by (\ref{first Thrall}), $W_3(2)\otimes_sW_3(2)\cong W_3(4)\oplus W_3(2^2)$.
Hence
\[
W_3(2^2)\otimes(W_3(2)\otimes_sW_3(2))\cong W_3(2^2)\otimes(W_3(4)\oplus W_3(2^2)),
\]
by (\ref{Young rule A})
\[
W_3(2^2)\otimes W_3(4)\cong W_3(6,2)
\oplus W_3(5,2,1)\oplus W_3(4,2^2),
\]
by (\ref{usage of LRrule})
\[
W_3(2^2)\otimes W_3(2^2)\cong W_3(4^2)\oplus W_3(4,3,1)\oplus W_3(4,2^2),
\]
and we obtain the decomposition for $W_3(2^2)\otimes(W_3(2)\otimes_sW_3(2))$.
In few cases we use also that
\[
W_d(\lambda_1+1,\ldots,\lambda_d+1)
\cong W_d(1^d)\otimes W_d(\lambda_1,\ldots,\lambda_d).
\]
For example,
\[
W_3(2,1^2)\otimes_sW_3(2,1^2)
\cong (W_3(1^3)\otimes W_3(1))\otimes_s(W_3(1^3)\otimes W_3(1))
\]
\[
=(W_3(1^3)\otimes W_3(1^3))\otimes(W_3(1)\otimes_sW_3(1))
\cong W_3(2^3)\otimes W_3(2)\cong W_3(4,2^2).
\]
\end{proof}

\begin{proposition}\label{hwv of degree 8}
For $d=3$, the following elements of $S=K[G_2\oplus\cdots\oplus G_6]$
are highest weight vectors:

\noindent For $\lambda=(4,3,1)$:
\[
w_1=\sum_{\sigma\in S_3}\text{\rm sign}(\sigma)
\text{\rm tr}([x_1,x_2]^2x_{\sigma(1)}x_{\sigma(2)})\text{\rm tr}(x_1x_{\sigma(3)}),
\]
\[
w_2=-\text{\rm tr}(s_3(x_1,x_2,x_3)x_1^2)\text{\rm tr}(x_1x_2^2)
\]
\[
+\text{\rm tr}(s_3(x_1,x_2,x_3)(x_1x_2+x_2x_1)\text{\rm tr}(x_1^2x_2)
-\text{\rm tr}(s_3(x_1,x_2,x_3)x_2^2)\text{\rm tr}(x_1^3),
\]
\[
w_3=\text{\rm tr}([x_1,x_2]^2)\text{\rm tr}(s_3(x_1,x_2,x_3)x_1),
\]
\[
w_4=\text{\rm tr}([x_1,x_2][x_1,x_3])(\text{\rm tr}(x_1^2)\text{\rm tr}(x_2^2)
\]
\[
-\text{\rm tr}^2(x_1x_2)) -\text{\rm tr}([x_1,x_2]^2)(\text{\rm tr}(x_1^2)\text{\rm tr}(x_2x_3)
-\text{\rm tr}(x_1x_2)\text{\rm tr}(x_1x_3)),
\]
\[
w_5=\text{\rm tr}(s_3(x_1,x_2,x_3)x_1)(\text{\rm tr}(x_1^2)\text{\rm tr}(x_2^2)-\text{\rm tr}^2(x_1x_2)),
\]
\[
w_6=(\text{\rm tr}(x_1^2x_2)\text{\rm tr}(x_2^2x_3)
-\text{\rm tr}(x_1^2x_3)\text{\rm tr}(x_2^3))\text{\rm tr}(x_1^2) +(-\text{\rm tr}(x_1^3)\text{\rm tr}(x_2^2x_3)
\]
\[
-\text{\rm tr}(x_1^2x_2)\text{\rm tr}(x_1(x_2x_3+x_3x_2))
+3\text{\rm tr}(x_1^2x_3)\text{\rm tr}(x_1x_2^2))\text{\rm tr}(x_1x_2)
\]
\[
+ (\text{\rm tr}(x_1^3)\text{\rm tr}(x_2^3)-\text{\rm tr}(x_1^2x_2)\text{\rm tr}(x_1x_2^2))\text{\rm tr}(x_1x_3)+
\]
\[
(\text{\rm tr}(x_1^3)\text{\rm tr}(x_1(x_2x_3+x_3x_2))
-2\text{\rm tr}(x_1^2x_2)\text{\rm tr}(x_1^2x_3))\text{\rm tr}(x_2^2)
\]
\[
+(-2\text{\rm tr}(x_1^3)\text{\rm tr}(x_1x_2^2) +2\text{\rm tr}^2(x_1^2x_2))\text{\rm tr}(x_2x_3),
\]
\[
w_7=\text{\rm tr}(s_3(x_1,x_2,x_3))(\text{\rm tr}(x_1^3)\text{\rm tr}(x_2^2)
-2\text{\rm tr}(x_1^2x_2)\text{\rm tr}(x_1x_2)+\text{\rm tr}(x_1x_2^2)\text{\rm tr}(x_1^2)).
\]

\noindent For $\lambda=(4,2^2)$:
\[
w_1=\text{\rm tr}(s_3(x_1,x_2,x_3)x_1^2)\text{\rm tr}(s_3(x_1,x_2,x_3)),
\]
\[
w_2=\sum_{\sigma\in S_3}\text{\rm sign}(\sigma)
\text{\rm tr}(s_3(x_1,x_2,x_3)x_{\sigma(1)}x_{\sigma(1)})\text{\rm tr}(x_1^2x_{\sigma(3)}),
\]
\[
w_3=\text{\rm tr}([x_1,x_2]^2)\text{\rm tr}([x_1,x_3]^2)-\text{\rm tr}^2([x_1,x_2][x_1,x_3]),
\]
\[
w_4=\text{\rm tr}^2(s_3(x_1,x_2,x_3)x_1),
\]
\[
w_5=\text{\rm tr}([x_2,x_3]^2)\text{\rm tr}^2(x_1^2)
-2\text{\rm tr}([x_1,x_3][x_2,x_3])\text{\rm tr}(x_1^2)\text{\rm tr}(x_1x_2)
\]
\[
+2\text{\rm tr}([x_1,x_2][x_2,x_3])\text{\rm tr}(x_1^2)\text{\rm tr}(x_1x_3)
+\text{\rm tr}([x_1,x_3]^2)\text{\rm tr}^2(x_1x_2)
\]
\[
-2\text{\rm tr}([x_1,x_2][x_1,x_3])\text{\rm tr}(x_1x_2)\text{\rm tr}(x_1x_3)
+\text{\rm tr}([x_1,x_2]^2)\text{\rm tr}^2(x_1x_3),
\]
\[
w_6=\text{\rm tr}([x_1,x_3]^2)(\text{\rm tr}(x_1^2)\text{\rm tr}(x_2^2)-\text{\rm tr}^2(x_1x_2))
\]
\[
-2\text{\rm tr}([x_1,x_2][x_1,x_3])*(\text{\rm tr}(x_1^2)\text{\rm tr}(x_2x_3)
-\text{\rm tr}(x_1x_2)\text{\rm tr}(x_1x_3))
\]
\[
+\text{\rm tr}([x_1,x_2]^2)(\text{\rm tr}(x_1^2)\text{\rm tr}(x_3^2)-\text{\rm tr}^2(x_1x_3)),
\]
\[
w_7=(-4\text{\rm tr}(x_1x_2^2)\text{\rm tr}(x_1x_3^2)+\text{\rm tr}^2(x_1(x_2x_3+x_3x_2)))\text{\rm tr}(x_1^2)
\]
\[
+ 4(2\text{\rm tr}(x_1^2x_2)\text{\rm tr}(x_1x_3^2)
-\text{\rm tr}(x_1^2x_3)\text{\rm tr}(x_1(x_2x_3+x_3x_2)))\text{\rm tr}(x_1x_2)
\]
\[
+4(2\text{\rm tr}(x_1^2x_3)\text{\rm tr}(x_1x_2^2)
-\text{\rm tr}(x_1^2x_2)\text{\rm tr}(x_1(x_2x_3+x_3x_2)))\text{\rm tr}(x_1x_3)
\]
\[
+4(-\text{\rm tr}(x_1^3)\text{\rm tr}(x_1x_3^2)+\text{\rm tr}^2(x_1^2x_3))\text{\rm tr}(x_2^2)
\]
\[
+4(\text{\rm tr}(x_1^3)\text{\rm tr}(x_1(x_2x_3+x_3x_2))
-2\text{\rm tr}(x_1^2x_2)\text{\rm tr}(x_1^2x_3))\text{\rm tr}(x_2x_3)
\]
\[
+4(\text{\rm tr}(x_1^3)\text{\rm tr}(x_1x_2^2)+\text{\rm tr}^2(x_1^2x_2))\text{\rm tr}(x_3^2),
\]
\[
w_8=\text{\rm tr}^2(s_3(x_1,x_2,x_3))\text{\rm tr}(x_1^2),
\]
\[
w_9=(\text{\rm tr}(x_1^2)(\text{\rm tr}(x_2^2)\text{\rm tr}(x_3^2)-\text{\rm tr}^2(x_2x_3))
\]
\[
-\text{\rm tr}^2(x_1x_2)\text{\rm tr}(x_3^2)+2\text{\rm tr}(x_1x_2)\text{\rm tr}(x_1x_3)\text{\rm tr}(x_2x_3)
-\text{\rm tr}^2(x_1x_3)\text{\rm tr}(x_2^2))\text{\rm tr}(x_1^2).
\]

\noindent For $\lambda=(3^2,2)$:
\[
w_1=-\text{\rm tr}([x_2,x_3]^2[x_1,x_2])\text{\rm tr}(x_1^2)
\]
\[
+\text{\rm tr}([x_1,x_2]([x_1,x_3][x_2,x_3]+[x_2,x_3][x_1,x_3]))\text{\rm tr}(x_1x_2)
\]
\[
-2\text{\rm tr}([x_1,x_2]^2[x_2,x_3])\text{\rm tr}(x_1x_3) -\text{\rm tr}([x_1,x_3]^2[x_1,x_2])\text{\rm tr}(x_2^2)
\]
\[
 +2\text{\rm tr}([x_1,x_2]^2[x_1,x_3])\text{\rm tr}(x_2x_3)-\text{\rm tr}([x_1,x_2]^3)\text{\rm tr}(x_3^2),
\]
\[
w_2=\text{\rm tr}(s_3(x_1,x_2,x_3)[x_1,x_2])\text{\rm tr}(s_3(x_1,x_2,x_3)),
\]
\[
w_3=\text{\rm tr}(s_3(x_1,x_2,x_3)x_1)\text{\rm tr}([x_1,x_2][x_2,x_3])
-\text{\rm tr}(s_3(x_1,x_2,x_3)x_2)\text{\rm tr}([x_1,x_2][x_1,x_3])
\]
\[ +\text{\rm tr}(s_3(x_1,x_2,x_3)x_3)\text{\rm tr}([x_1,x_2]^2),
\]
\[
w_4=\text{\rm tr}(s_3(x_1,x_2,x_3)x_1)(\text{\rm tr}(x_1x_2)\text{\rm tr}(x_2x_3)
-\text{\rm tr}(x_1x_3)\text{\rm tr}(x_2^2))
\]
\[
+\text{\rm tr}(s_3(x_1,x_2,x_3)x_2)(-\text{\rm tr}(x_1^2)\text{\rm tr}(x_2x_3)
+\text{\rm tr}(x_1x_2)\text{\rm tr}(x_1x_3))
\]
\[
+\text{\rm tr}(s_3(x_1,x_2,x_3)x_3)(\text{\rm tr}(x_1^2)\text{\rm tr}(x_2^2)-\text{\rm tr}^2(x_1x_2)).
\]
In each of the cases,
every highest weight vector $w\in W_3(\lambda)\subset \omega^2(S)$ is equal
to a linear combination of $w_i$.
\end{proposition}

\begin{proof}
The considerations use the decomposition of
$(\omega^2(S))^{(8)}$ given in Proposition \ref{decomposition of degree 8}.
The highest weight vectors have been found as those in Propositions
\ref{hwv of degree up to 6} and \ref{hwv of degree 7}. Of course, if
we already know the explicit form of the (candidates for) highest weight vectors $w_i$,
we can check that they are linearly independent in $\omega^2(S)$
and satisfy the requirements of Lemma \ref{criterion for hwv}. Since for each
$\lambda$ the number of the highest weight vectors $w_i$ coincides with the
multiplicity of $W_3(\lambda)\subset \omega^2(S)$ from
Proposition \ref{decomposition of degree 8}, we conclude that every
highest weight vector $w\in W_3(\lambda)\subset \omega^2(S)$ is equal
to a linear combination of $w_i$.
\end{proof}

Finally, we shall calculate the Hilbert series of the kernel
of the natural homomorphism $S\to C_0$ for $d=3$.

\begin{lemma}\label{Hilbert series of kernel}
Let $d=3$ and let $J$ be the kernel of the natural homomorphism
$S\to C_0$. Let the Hilbert series of $J$ be
\[
H(J,t_1,t_2,t_3)=\sum_{k\geq 0}h_k(t_1,t_2,t_3),
\]
where $h_k$ is the homogeneous component of degree $k$ of $H(J,t_1,t_2,t_3)$. Then
$h_k=0$ for $k\leq 6$,
\[
h_7=S_{(3,2^2)}(t_1,t_2,t_3),
\]
\[
h_8=S_{(4,3,1)}(t_1,t_2,t_3) + 2S_{(4,2^2)}(t_1,t_2,t_3) + S_{(3^2,2)}(t_1,t_2,t_3).
\]
\end{lemma}

\begin{proof}
Clearly, the Hilbert series of the kernel $J$ is equal to the difference
of the Hilbert series of $S$ and $C_0$. For $d=3$  we have
that
\[
G_2\oplus\cdots\oplus G_6=W(2)\oplus W(3)\oplus W(1^3)
\]
\[
\oplus W(2^2)\oplus W(2,1^2)
\oplus W(3,1^2)\oplus W(2^2,1)\oplus W(3^2)
\]
and its Hilbert series is
\[
H(G_2\oplus\cdots\oplus G_6,t_1,t_2,t_3)=
\sum a_{k_1k_2k_3}t_1^{k_1}t_2^{k_2}t_3^{k_3}
\]
\[
=S_{(2)}+S_{(3)}+S_{(1^3)}+S_{(2^2)}+S_{(2,1^2)}
+S_{(3,1^2)}+S_{(2^2,1)}+S_{(3^2)}.
\]
The Hilbert series of the symmetric algebra is
\[
H(S,t_1,t_2,t_3)=\prod\frac{1}{(1-t_1^{k_1}t_2^{k_2}t_3^{k_3})^{a_{k_1k_2k_3}}}.
\]
Now the result follows by evaluation of the coefficients $a_{k_1k_2k_3}$ and
expanding the first several homogeneous components of the difference of the
Hilbert series of $S$ and of the Hilbert series of $C_0$, which is given in (\ref{Hilbert series of C0}).
\end{proof}

\begin{corollary}\label{module of relations of degree 7 and 8}
For $d=3$, the algebra $C_0$ has a minimal system of defining relations
with the property that the relations of degree $7$ and $8$ form $GL_3$-modules
isomorphic, respectively, to $W_3(3,2^2)$ and $W_3(4,3,1) + 2W_3(4,2^2) + W_3(3^2,2)$.
\end{corollary}

\begin{proof}
If $J$ is the kernel of the natural homomorphism $S\to C_0$, then a minimal homogeneous system of generators
of $J$ is obtained as a factor space of $J$ modulo $J\omega(S)$. Since $\omega(S)$ contains no homogeneous
elements of degree 1, we obtain that the multihomogeneous components of total degree 7 and 8 of $J$ and $J/J\omega(S)$
are of the same dimension. Hence $J^{(7)}$ and $J^{(8)}$ are isomorphic as $GL_3$-modules
to $(J/J\omega(S))^{(7)}$ and $(J/J\omega(S))^{(8)}$, respectively, and the conclusion follows from
the expressions of $h_7$ and $h_8$ given in Lemma \ref{Hilbert series of kernel}.
\end{proof}

\section{Main results}

Now we present the explicit defining relations of degree 7 of the algebra $C_{3d}$ for any $d\geq 3$
and of degree 8 for the algebra $C_{33}$, with respect to the generators of Abeasis and Pittaluga \cite{AP}.
As we already mentioned, by (\ref{replacing with traceless matrices})
it is sufficient to give the defining relations of the algebra $C_0$ generated by
traces $\text{tr}(x_{i_1}\cdots x_{i_k})$ of products of the traceless matrices $x_i$.
As in the previous sections, we denote by $S$ the symmetric algebra of the $GL_d$-module
$G_2\oplus\cdots\oplus G_6$ of generators of $C_0$ and call defining relations of $C_0$
the expressions $f=0$, where $f$ is an element of the kernel $J$ of
the natural homomorphisms $S\to C_0$.

\begin{theorem}\label{relations of degree 7}
Let $d\geq 3$. The algebra $C_0$ does not have any defining relations of degree $\leq 6$.
The $GL_d$-module structure of the homogeneous defining relations of degree $7$
of $C_0$, i.e., of the component $J^{(7)}$ in $S$ is
\[
J^{(7)}=W_d(4,1^3)\oplus W_d(3,2^2)\oplus W_d(3,2,1^2)\oplus W_d(2^3,1)
\oplus W_d(2^2,1^3)\oplus W_d(2,1^5).
\]
In the notation of Proposition \ref{hwv of degree 7},
the defining relations of $C_0$ which are highest weight vectors are:

\noindent For $\lambda=(4,1^3)$:
\begin{equation}\label{relation hwv 4111}
12w_1-15w_2-20w_3=0;
\end{equation}

\noindent For $\lambda=(3,2^2)$:
\[
2w_1-w_2+2w_3=0.
\]

\noindent For $\lambda=(3,2,1^2)$:
\[
-6w_1+10w_3-15w_4+40w_6=0.
\]

\noindent For $\lambda=(2^3,1)$:
\[
12w_1+w_2=0.
\]

\noindent For $\lambda=(2^2,1^3)$:
\[
w_2=0.
\]

\noindent For $\lambda=(2,1^5)$:
\[
2w_1-5w_2=0.
\]
\end{theorem}

\begin{proof}
In all cases the idea is the same. We already know that there are no relations of
degree $\leq 11$ for $d=2$ and the only $GL_3$-module of relations is isomorphic to $W_3(3,2^2)$.
Hence, we have to consider the cases in Proposition \ref{hwv of degree 7} only.

We shall consider in detail the case $\lambda=(4,1^3)$. The possible relations $w=0$
are linear combinations of $w_1,w_2,w_3$. We assume that
\begin{equation}\label{relation 4111}
w=\xi_1w_1+\xi_2w_2+\xi_3w_3=0
\end{equation}
and evaluate $w$ on the traceless matrices (\ref{first matrix}) and (\ref{other matrices}).
The coefficients of the monomials $(x_{11}^{(1)})^4x_{12}^{(2)}x_{22}^{(3)}x_{21}^{(4)}$ and
$(x_{11}^{(1)})^4x_{13}^{(2)}x_{32}^{(3)}x_{21}^{(4)}$ are, respectively,
$20\xi_1-8\xi_2+18\xi_3$ and $20\xi_1+12\xi_3$. Hence the equation (\ref{relation 4111}) implies that
\[
20\xi_1-8\xi_2+18\xi_3=20\xi_1+12\xi_3=0.
\]
Up to a multiplicative constant, the only solution of this system is
\[
\xi_1=12,\quad \xi_2=-15,\quad \xi_3=-20.
\]
Hence, there is only one possible candidate for a defining relation which is a highest weight vector
of some $W_d(4,1^3)$. We evaluate once again (\ref{relation 4111})
on (\ref{first matrix}) and (\ref{other matrices})
for these values of $\xi_1,\xi_2,\xi_3$ and obtain that $w(x_1,x_2,x_3,x_4)=0$. Hence the multiplicity
of $W_d(4,1^3)$ in $J$ is equal to 1 and the corresponding relation is (\ref{relation hwv 4111}).
We want to mention that the case $\lambda=(3,1^4)$ does not participate in the statement
of the theorem, because the multiplicity of $W_d(3,1^4)$ in $J$ is 0.
\end{proof}

\begin{corollary}
The dimension of the defining relations of degree $7$ of the algebra $C_{3d}$ is equal to
\[
r_7=r_7(d)=\frac{2}{7!}(d+1)d(d-1)(d-2)(41d^3-86d^2+114d-360).
\]
\end{corollary}

\begin{proof}
Since $C_{3d}$ does not satisfy relations of degree $\leq 6$,
and the constants are the only elements of degree 0 in
$K[\text{tr}(X_1),\ldots,\text{tr}(X_d)]$, the dimension of the relations of degree 7
of $C_{3d}$ coincides with this dimension in $C_0$. Now the proof is complete using
Theorem \ref{relations of degree 7} and the dimension formula
(\ref{dim of W}) for $W_d(\lambda)$.
\end{proof}

\begin{theorem}\label{relations of degree 8}
Let $d=3$. The $GL_d$-module structure of the homogeneous component $J^{(8)}$ of degree $8$ in $S$ is
\[
J^{(8)}=W_3(4,3,1)\oplus 2W_3(4,2^2)\oplus W_3(3,2^2,2).
\]
In the notation of Proposition \ref{hwv of degree 8},
the defining relations which are highest weight vectors are:

\noindent For $\lambda=(4,3,1)$:
\begin{equation}\label{relation hwv 431}
-6w_1-18w_2+3w_3+3w_5-8w_7=0;
\end{equation}

\noindent For $\lambda=(4,2^2)$: All nontrivial linear combinations of
\[
w_1-15w_2+3w_3+\frac{21}{4}w_4-\frac{5}{2}w_5+\frac{5}{2}w_6-3w_7+2w_9=0,
\]
\[
-36w_2+6w_3+\frac{27}{2}w_4-6w_5+6w_6-9w_7+w_8+6w_9=0.
\]

\noindent For $\lambda=(3^2,2)$:
\[
6w_1+2w_2-3w_3-3w_4=0.
\]
The number $r_8$ of the defining relations of degree $8$ of any homogeneous minimal system
of defining relations of the algebra $C_{33}$ is equal to $30$.
\end{theorem}

\begin{proof}
The decomposition of $J^{(8)}$ is given in
Corollary \ref{module of relations of degree 7 and 8}. The explicit form of the highest weight
vectors is obtained as in the proof of Theorem \ref{relations of degree 7}. The number of
defining relations of degree 8 in any homogeneous minimal system of defining relations for $C_{33}$
is equal to the dimension of the relations $J^{(8)}$ of $C_0$. For the proof that $r_8=30$
it is sufficient to use the dimension formula (\ref{dim of W}) for $W_d(\lambda)$.
\end{proof}

\begin{remark}
Using Lemma \ref{basis of module} we can find an explicit basis of the set of defining relations
of degree 7 for $C_0$, $d\geq 3$, and of degree 8 for $C_0$, $d=3$.
\end{remark}

\section*{Acknowledgements}

This project was started when the second author visited the
University of Pa\-ler\-mo.
He is very grateful for the hospitality and the creative
atmosphere during his stay there.

\end{document}